\documentclass[authoryear,final,5p,times,twocolumn]{elsarticle}

\usepackage{amssymb}
\usepackage{siunitx}
\usepackage[colorlinks = true,
linkcolor = blue,
urlcolor  = blue,
citecolor = blue,
anchorcolor = blue]{hyperref}
\usepackage{physics} 
\usepackage{bm} 
\usepackage{stfloats} 
\usepackage[capitalise]{cleveref}
\usepackage{tikz}
\usepackage{svg}
\usetikzlibrary{arrows,arrows.meta,calc,decorations.markings,decorations.markings,math,patterns,positioning,shadows.blur,shapes,shapes.symbols,spy}
\usepackage{subcaption}
\usepackage[colorinlistoftodos]{todonotes}
\usepackage{xurl}
\usepackage{framed}

\journal{Applied Thermal Engineering}

\begin{document}
	\let\WriteBookmarks\relax
	\def\floatpagepagefraction{1}
	\def\textpagefraction{.001}
	
	\newcommand{\NurbsBasis}{\hat{N}}
	\newcommand{\BSplineBasis}{\hat{B}}
	\newcommand{\NurbsDegree}{p}
	\newcommand{\NumberBasisFunctions}{n}
	\newcommand{\NurbsWeight}{w}
	\newcommand\st{\mathrm{st}}
	\newcommand\rt{\mathrm{rt}}
	\newcommand\ag{\mathrm{ag}}
	\newcommand\IntegG{\mathrm{d}\Gamma}
	\newcommand\Integ{\mathrm{d}\Omega}
	\newcommand{\opt}{\mathrm{opt}}
	\newcommand{\Bhigh}{B_{\mathrm{high}}}
	\newcommand{\Blow}{B_{\mathrm{low}}}
	
	\newcommand{\OmegaRotor}{\Omega_\rt}
	\newcommand{\OmegaStator}{\Omega_\st}
	\newcommand{\GammaAirGap}{\Gamma_\ag}
	
	\newcommand{\NAMR}{n_\mathrm{AMR}}
	\newcommand{\rAMR}{r_\mathrm{AMR}}
	
	\definecolor{TUDa-2a}{HTML}{009CDA}
	\definecolor{TUDa-2b}{HTML}{0083CC}
	\definecolor{TUDa-3d}{HTML}{0071F3}
	\definecolor{TUDa-3a}{HTML}{50B695}
	\definecolor{TUDa-9b}{HTML}{E6001A}

\begin{frontmatter}



\title{How to Build the Optimal Magnet Assembly for Magnetocaloric Cooling: Structural Optimization with Isogeometric Analysis}


\author[1]{Michael Wiesheu}

\ead{michael.wiesheu@tu-darmstadt.de}

\author[1]{Melina Merkel}
\author[2]{Tim Sittig}
\author[2]{Dimitri Benke}
\author[2]{Max Fries}
\author[1]{Sebastian Schöps}
\author[3]{Oliver Weeger}
\author[4]{Idoia {Cortes Garcia}}

\affiliation[1]{organization={Computational Electromagnetics Group, Technical University of Darmstadt},
	addressline={Schloßgartenstr.~8}, 
	city={Darmstadt},
	postcode={64289}, 
	country={Germany}}

\affiliation[2]{organization={MagnoTherm Solutions GmbH},
	addressline={Pfungstädter Str.~102}, 
	city={Darmstadt},
	postcode={64297}, 
	country={Germany}}

\affiliation[3]{organization={Cyber-Physical Simulation, Technical University of Darmstadt},
	addressline={Dolivostr.~15}, 
	city={Darmstadt},
	postcode={64293}, 
	country={Germany}}

\affiliation[4]{organization={Department of Mechanical Engineering, Eindhoven University of Technology},
	addressline={Gemini-Zuid 0.143 PO Box 513}, 
	city={Eindhoven},
	postcode={5600 MB}, 
	country={The Netherlands}}

\begin{abstract}
	In  the  search  for  more  efficient  and  less environmentally harmful cooling technologies, the  field  of  magnetocalorics  is  considered  a  promising  alternative. To generate cooling spans, rotating permanent magnet assemblies are used to  cyclically  magnetize  and  demagnetize  magnetocaloric  materials,  which  change  their  temperature  under  the  application of  a  magnetic  field.  In  this  work,  an  axial  rotary  permanent  magnet assembly,  aimed  for commercialization, is computationally designed using topology and shape optimization. This is efficiently facilitated in an isogeometric analysis framework, where harmonic mortaring is applied to couple the rotating rotor-stator system of the multipatch model.  Inner, outer and co-rotating assemblies are compared and optimized designs for different magnet masses are determined.  These simulations are used to homogenize the magnetic flux density in the magnetocaloric material. The resulting torque is analyzed for different geometric parameters. Additionally, the influence of anisotropy in the active magnetic regenerators is studied in order to guide the magnetic flux. Different examples are analyzed and classified to find an optimal magnet assembly for magnetocaloric cooling.
	
\end{abstract}


\end{frontmatter}

\section{Introduction}
As of 2018, cooling contributed to \SI{17}{\percent} of the worldwide electricity consumption and this demand is projected to triple by 2050 \citep{UnitedNationsEntironmentProgramme2020}.
Furthermore,  hydrofluorocarbons (HFCs) -- highly potent greenhouse gases -- are commonly used as refrigerants in conventional vapor-compression cooling systems. Combining direct emissions caused by the leakage of HFC refrigerant gases and indirect emissions from the power generation, cooling accounts for \SI{8}{\percent} of the global CO$_2$ equivalent emissions\footnote{Calculated with data from \citep{GHGdata} and \citep{GreenCoolingInitiative} in 2016.}. 
Switching to climate-friendly refrigerants and improving the energy efficiency of cooling appliances are the key challenges the cooling sector is facing in the fight against climate change \citep{UnitedNationsEntironmentProgramme2020}.
Magnetocaloric cooling involves no environmentally critical gases and has the potential to increase efficiency
, as supported by measurements from \cite{PerformanceOfMagneticRefrigerator}.
Temperature spans are generated with magnetocaloric materials (MCMs), which change their temperature when a magnetic field is applied or removed. This is commonly achieved  by rotating permanent magnet assemblies, that cyclically magnetize and demagnetize the MCM. The heat transfer is accomplished by a heat transfer fluid -- usually water -- which is pumped according to the magnetization and demagnetization phases.

A variety of different permanent magnet assemblies are known from the literature \citep{BJORK2010437,MagnetocaloricEnergyConversion}. 
Such assemblies are commonly based on Halbach (cylinder) arrays \citep{Aprea2018,Arnold2014,Bjrk2008} or optimized magnet systems \citep{optimizedrefrigeration,PermanentMagnetDesignForMagneticHeatPumps}. 

In this work, the interest lies particularly in the structural optimization of an axial rotary system design, as initially proposed by \cite{Okamura}. Axial rotary systems have the advantageous property of linear scalability in the axial direction (neglecting edge effects), making them a favorable choice for commercialization.

At present, commercial vapor-compression systems are more than a hundred years ahead in research and development \citep{Botoc2021}.
Thus, simulation and optimization are fundamental tools to rapidly improve magnetocaloric cooling devices and make them competitive in terms of power density and efficiency. 
Here, this is achieved with the use of Isogeometric Analysis (IGA). Being able to model geometries exactly and to apply shape optimization in a straightforward way are particular advantages of IGA. As a result, IGA has successfully been applied to model and optimize rotating machines \citep{RecentAdvancesOfIsogeometricAnalysis}.

The structure of this paper is as follows: \cref{chap:methodology} gives insights to the underlying geometrical and physical modeling. In \cref{chap:StructuralOptimization}, an optimized magnet design for the rotor is generated with topology and shape optimization. \cref{chap:FinalDesignStudies} presents studies for the torque and the magnetic flux density. The work is concluded in \cref{chap:Conclusion}.

\begin{table*}[!ht]   
\begin{framed}
\textbf{Nomenclature} \\

\begin{minipage}{0.5\linewidth}
    \textit{Abbreviations}\\
     \begin{tabular}{l l}
         AMR & Active Magnetic Regenerator \\
         FEM & Finite Element Method \\
         HFC & Hydrofluorocarbon \\
         IGA & Isogeometric Analysis \\
         MCM & Magnetocaloric Material \\
         NURBS & Non-Uniform Rational B-Spline \\
         PDE & Partial Differential Equation \\
         SIMP & Solid Isotropic Material with Penalization \\
         La-Fe-Si & Alloy based on Lanthanum, Iron and Silicium
     \end{tabular}\\
     
     \textit{Roman symbols} \\
     \begin{tabular}{l l}
         $\mathbf{A}$ & Magnetic vector potential (\SI{}{V.s.m^{-1}})\\
         A & Area (\SI{}{\meter\squared})\\
         $\mathbf{B}$ & Magnetic flux density (\SI{}{T}) \\
         $\mathbf{b}$ & Discretized remanence\\
         $\mathbf{F}$ & Mapping \\
         $f$ & Phase fraction \\
         $\mathbf{G}$ & Coupling Matrix \\
         $\mathbf{H}$ & Magnetic field strength (\SI{}{A.m^{-1}}) \\
         $\mathbf{J}$ & Electric current density (\SI{}{A.m^{-2}}) \\
         $\mathbf{K}$ & Stiffness Matrix \\
         $\mathbf{n}$ & Surface normal vector \\
         $N$ & Basis function \\
         $p$ & Polynomial degree \\
         $\mathbf{P}$ & Control point \\
         $q$ & Penalization factor \\
         $\mathbf{R}$ & Rotation matrix \\
         $\mathbf{u}$ & Discretized magnetic vector potential\\
         $V$ & Volume \\
         $w$ & NURBS weight \\
     \end{tabular}

\end{minipage}
\begin{minipage}{0.5\linewidth}

    \textit{Greek symbols}\\
     \begin{tabular}{l l}
         $\alpha$ & Phase angle \\
         $\beta$ & Rotation angle \\
         $\Gamma$ & Domain surface \\
         $\bm{\lambda}$ & Lagrangian multiplier \\
         $\bm{\mu}$ & Magnetic permeability (\SI{}{H.m^{-1}}) \\
         $\bm{\nu}$ & Magnetic reluctivity (\SI{}{m.H^{-1}}) \\
         $\Omega$ &  Physical domain \\
         $\hat{\Omega}$ & Parametric domain \\
         $\theta$ & Evaluation angle in the air gap \\
         $\Xi$ & Knot vector \\
         $\xi$ & Knot vector entry \\
     \end{tabular}\\
    
    \textit{Subscripts}\\
     \begin{tabular}{l l}
         ag & air gap \\
         azi & azimuthal \\
         eval & evaluation \\
         high & high field \\
         low & low field \\
         mag & magnetic \\
         max & maximum \\
         opt & optimization \\
         $\%$ & Percentage share \\
         pm & permanent magnet \\
         r & remanent \\
         rad & radial \\
         rel & relative \\
         rt & rotor \\
         st & stator \\
     \end{tabular}\\
     
\end{minipage}


\end{framed}
\end{table*}

\section{Methodology} 
\label{chap:methodology}

\subsection{Geometric modeling}\label{sec:geometry}
For the simulation of physical problems, the geometry needs to be represented. Classically, in the Finite Element Method (FEM), this is often done using a mesh of lowest order, e.g., triangles or tetrahedra. More detailed geometries, i.e., a finer mesh, lead to higher computational efforts and less detailed ones to larger geometrical and numerical errors. Thus, there is always a trade-off between computational cost and model accuracy. 
One way to resolve this problem is to use the CAD basis functions, i.e., B-Splines and Non-Uniform Rational B-Splines (NURBS), to represent complex geometrical structures exactly even within the FEM.

B-Spline basis functions $\BSplineBasis_i^p$ are defined by their degree $\NurbsDegree$ as well as the knot vector 
\begin{equation}
    \Xi =\big\{ \xi _{1} ,\ \xi _{2} ,\ ...,\ \xi _{n+p+1} \big\},
\end{equation}
and can be calculated using Cox-de Boor's recursion formula~\citep{TheNurbsBook}. The number of basis functions $\NumberBasisFunctions$ results from $\NurbsDegree$ and $\Xi$.
NURBS basis functions are constructed by weighting each B-Spline basis function $\BSplineBasis_i^p$ with a additional weight $\NurbsWeight_i$ as 
\begin{equation}
    \NurbsBasis_{i}^{\NurbsDegree}( \hat{x} ) \ =\ \frac{\NurbsWeight_{i} \BSplineBasis_{i}^{\NurbsDegree}( \hat{x} )}{\sum_{j=1}^\NumberBasisFunctions \NurbsWeight_{j} \BSplineBasis_{j}^{\NurbsDegree}( \hat{x} )}, \label{eq:NumericalMethods:Nurbs:NURBS}
\end{equation}
for $\hat{x}\in\hat{\Omega}\subset\mathbb{R}$, where $\hat{\Omega}=[\xi _{1}, \xi _{n+p+1}]$ is the reference domain.
Multivariate B-Spline and NURBS basis functions can be defined using a tensor-product approach.  
To describe a 1-dimensional manifold $\Omega\subset\mathbb{R}^d$ in the $d$-dimensional physical space  as a NURBS curve, each basis function is associated with a control point $\mathbf{P}_i\in\mathbb{R}^d$ via multiplication. The curve parametrization is then given by
\begin{equation}
    \mathbf{F}:\hat\Omega\to\Omega, \quad
    \mathbf{x} = \mathbf{F}( {\hat{x}} ) :=\sum _{i=1}^{\NumberBasisFunctions}\mathbf{P}_{i} \NurbsBasis_{i}^{\NurbsDegree}( {\hat{x}} ). 
    \label{eq:NumericalMethods:Nurbs:3DCurve}
\end{equation}
Surfaces and volumes can be defined similarly using multivariate basis functions and a mesh of control points.

\begin{figure}[h]
    \centering
    \includegraphics[trim= 150 40 200 15, clip, width=\linewidth]{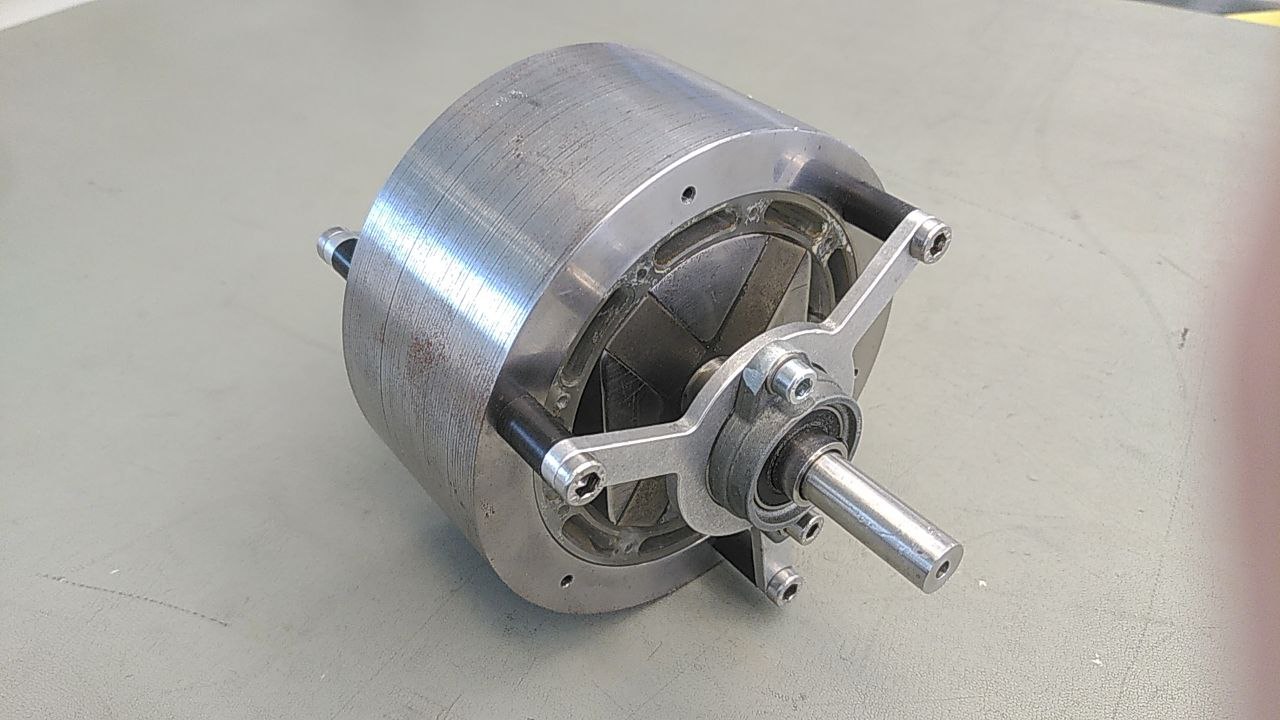}
    \caption{Exemplary magnet assembly of an axial rotary cooling device. Published with kind permission from Christian Vogel. }
    \label{fig:Methodology:PhysicalDevice}
\end{figure}

This allows for the exact representation of conic sections, such as circles or ellipses, which is particularly useful for rotating machines, such as the one shown in \cref{fig:Methodology:PhysicalDevice}. 
The working principle of this exemplary axial rotary cooling device will be explained in the following chapters.

\subsection{Magnetostatics} 
When modeling magnetic machines, it is often sufficient to only consider simplified approximations of Maxwell's equations.
Here, nonlinear material behavior and time dependent phenomena such as eddy currents are neglected.
The magnetic vector potential formulation \citep{Jackson_1998aa}
\begin{align}
    \nabla \times \left( \bm{\nu}( \nabla \times \mathbf{A} -\mathbf{B}_{\mathrm{r}})\right) &= \mathbf{J},
    \label{eq:methodology:StrongMagnetostatic}
\end{align}
then describes the magnetic field, where the magnetic vector potential $\mathbf{A}$ is defined such that the magnetic flux density is
\begin{align}
    \mathbf{B}&=\nabla\times\mathbf{A}. \label{eq:methodology:MagVecPot}
\end{align}
Furthermore, $\bm{\nu} = \bm{\mu}^{-1}$ is the (potentially anisotropic) reluctivity, i.e., the inverse permeability, $\mathbf{B}_{\mathrm{r}}$ the remanence of the permanent magnets and $\mathbf{J}$ the electric current density. Here, all variables are assumed to be space-dependent. The linear material law
\begin{equation}
    \mathbf{B} \ =\ \bm{\mu}\mathbf{H} \ +\ \mathbf{B}_{\mathrm{r}} \label{eq:methodology:constitutiveLaw}
\end{equation}
is used to model the relation between the magnetic field strength $\mathbf{H}$ and magnetic flux density $\mathbf{B}$ considering the remanence of the material.

In the following, we assume that the geometry is sufficiently long in axial direction and has no cross-sectional changes, such that magnetic fields in the $z$-direction are negligible. We therefore restrict ourselves to two dimensional magnetostatic modeling and simulations.
\subsubsection{Magnetostatics for 2D}
By assuming
$\mathbf{A} = \begin{bmatrix}
0 & 0 & A_{z}
\end{bmatrix}^{\mathrm{\top}}$ for the 2D-case,
(\ref{eq:methodology:StrongMagnetostatic}) simplifies to the Poisson problem \citep{Salon_1995aa}
\begin{equation}
    -\nabla \cdotp \left(
    \bm{\nu}^\bot\left( \nabla A_{z} -\mathbf{B}_{\mathrm{r}}^{\bot }\right)\right) =J_{z}
    \label{eq:methodology:strong2DAnisotropic}
\end{equation}
with the anisotropic reluctivity 
\begin{align}
    \bm{\nu}^\bot= \frac{\bm{\mu} ^{\top}}{\det(\bm{\mu})}\in \mathbb{R}^{2\times 2}, \label{eq:reluctivity2D}
\end{align}
where $\bm{\mu}$ is also reduced to 2D. In the isotropic case, (\ref{eq:reluctivity2D}) further simplifies to $\nu=\mu^{-1}\in\mathbb{R}$.
The variable $\mathbf{B}_{\mathrm{r}}^{\bot }$ contains the  perpendicular components of the remanence and is defined as
\begin{equation}
    \mathbf{B}_{\mathrm{r}}^{\bot } =\begin{pmatrix}
    -B_{\mathrm{r}y}\\
    B_{\mathrm{r}x}
    \end{pmatrix}. \label{eq:methodology:BremPerpendicular}
\end{equation}
The magnetic flux density is then calculated from (\ref{eq:methodology:MagVecPot}) as
\begin{equation}
\mathbf{B}=\begin{pmatrix}
\partial A_{z} /\partial y\\
-\partial A_{z} /\partial x
\end{pmatrix}.
\end{equation}

\subsubsection{Magnetocaloric material model}

Implementing an  active magnetic regenerator (AMR) is the most common way to create larger temperature spans in magnetocaloric cooling. The AMR consists of a small-scale MCM structure which is in contact with the heat transfer fluid. The geometry of the MCM is either classified as \textit{packed bed} or \textit{ordered structure}  \citep{EnergyApplicationsOfMagnetocaloricMaterials}. Common examples are packed spheres or parallel plates, respectively.
The feature size of the MCM is usually below \SI{1}{\mm} in order to ensure a good heat transfer between the solid material and the perfusing fluid. 

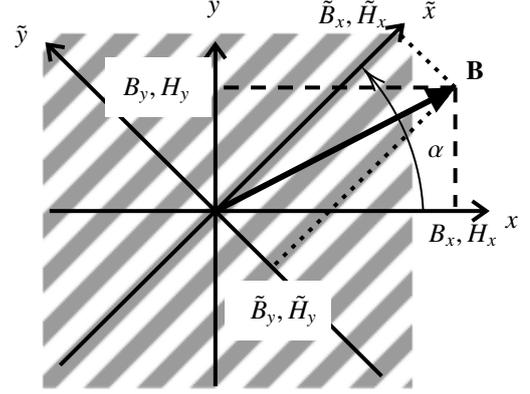
\begin{figure}
    \centering    
        \tikzset{every picture/.style={line width=0.75pt}} 
        
        \begin{tikzpicture}[x=0.75pt,y=0.75pt,yscale=-1,xscale=1]
        
        \draw (336.5,141.5) node  {\includegraphics[width=147.75pt,height=147.75pt]{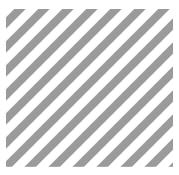}};
        \draw [line width=1.5]  (247.47,142.1) -- (465.47,142.1)(330.24,58.4) -- (330.24,230.2) (458.47,137.1) -- (465.47,142.1) -- (458.47,147.1) (325.24,65.4) -- (330.24,58.4) -- (335.24,65.4)  ;
        \draw [line width=2.25]    (330.24,142.1) -- (445.56,82.3) ;
        \draw [shift={(450,80)}, rotate = 152.59] [fill={rgb, 255:red, 0; green, 0; blue, 0 }  ][line width=0.08]  [draw opacity=0] (14.29,-6.86) -- (0,0) -- (14.29,6.86) -- cycle    ;
        \draw [line width=1.5]  [dash pattern={on 5.63pt off 4.5pt}]  (450,80) -- (450,140) ;
        \draw [line width=1.5]  (252.74,219.6) -- (423.63,48.7)(246.54,58.4) -- (413.27,225.14) (415.15,50.12) -- (423.63,48.7) -- (422.22,57.19) (247.95,66.89) -- (246.54,58.4) -- (255.02,59.81)  ;
        \draw [line width=1.5]  [dash pattern={on 5.63pt off 4.5pt}]  (450,80) -- (330,80) ;
        \draw [line width=1.5]  [dash pattern={on 1.69pt off 2.76pt}]  (450,80) -- (421,52) ;
        \draw [line width=1.5]  [dash pattern={on 1.69pt off 2.76pt}]  (360,169) -- (450,80) ;
        \draw    (434,143) .. controls (433.02,117.52) and (423.4,92.51) .. (406.07,72.69) ;
        \draw [shift={(405,71.48)}, rotate = 48.01] [color={rgb, 255:red, 0; green, 0; blue, 0 }  ][line width=0.75]    (10.93,-3.29) .. controls (6.95,-1.4) and (3.31,-0.3) .. (0,0) .. controls (3.31,0.3) and (6.95,1.4) .. (10.93,3.29)   ;
        \draw  [draw opacity=0][fill={rgb, 255:red, 255; green, 255; blue, 255 }  ,fill opacity=1 ] (275,65) -- (328,65) -- (328,95) -- (275,95) -- cycle ;
        \draw  [draw opacity=0][fill={rgb, 255:red, 255; green, 255; blue, 255 }  ,fill opacity=1 ] (335,177) -- (388,177) -- (388,207) -- (335,207) -- cycle ;
        
        \draw (473,143.4) node [anchor=north west][inner sep=0.75pt]    {$x$};
        \draw (325,35) node [anchor=north west][inner sep=0.75pt]    {$y$};
        \draw (432,35.4) node [anchor=north west][inner sep=0.75pt]    {$\tilde{x}$};
        \draw (229,47.4) node [anchor=north west][inner sep=0.75pt]    {$\tilde{y}$};
        \draw (435,147.4) node [anchor=north west][inner sep=0.75pt]    {$B_{x} ,H_{x}$};
        \draw (282.5,73.4) node [anchor=north west][inner sep=0.75pt]    {$B_{y} ,H_{y}$};
        \draw (380,35) node [anchor=north west][inner sep=0.75pt]    {$\tilde{B}_{x} ,\tilde{H}_{x}$};
        \draw (346,183.4) node [anchor=north west][inner sep=0.75pt]    {$\tilde{B}_{y} ,\tilde{H}_{y}$};
        \draw (435,107.4) node [anchor=north west][inner sep=0.75pt]    {$\alpha $};
        \draw (454,65.4) node [anchor=north west][inner sep=0.75pt]    {$\mathbf{B}$};
        
        \end{tikzpicture}

    \caption{Visualization of a rotated two phase material in an arbitrary magnetic field. 
    }
    \label{fig:modeling:multiscale:rotated}
\end{figure}

To account for the anisotropic behavior of such microstructures, the analytic Voigt and Reuss bounds for the permeability are used to determine the edge cases of the material behavior \citep{Salon_1995aa}. The Voigt-case corresponds to a layered material with the different phases arranged in the direction of the magnetic field and defines the upper limit for the permeability as
\begin{equation}
    \mu_{\mathrm{Voigt}} =( f_{1} \mu _{1} +f_{2} \mu _{2}) = \left< \mu \right>, \label{eq:modeling:multiscale:voigt2}
\end{equation}
where $f_1$ and $f_2$ are the volume fractions and $\left< \cdot \right>$ denotes the averaged value. 
The Reuss-case corresponds to a layered material with the phases arranged perpendicular to the magnetic field and defines the lower bound for the permeability as
\begin{equation}
    \mu_{\mathrm{Reuss}} \ =\ \left(\frac{f_{1}}{\mu _{1}} +\frac{f_{2}}{\mu _{2}}\right)^{-1} =\left< \mu ^{-1}\right> ^{-1}. \label{eq:modeling:multiscale:reuss2}
\end{equation}
In the case of parallel phases rotated by the angle $\alpha$, as shown in Fig. \ref{fig:modeling:multiscale:rotated}, the effective permeability tensor $\bm{\mu}$ is obtained by applying
\begin{equation}
    \bm{\mu } =\mathbf{R}_{\alpha }\tilde{\bm{\mu }}\mathbf{R}_{\alpha }^{-1} \label{eq:modeling:MuTensor}.
\end{equation}
This describes the coordinate transformation to the known unrotated reference system by the rotation matrix
\begin{equation}
    \mathbf{R}_{\alpha} = \begin{pmatrix}
    \cos( \alpha ) & -\sin( \alpha )\\
    \sin( \alpha ) & \cos( \alpha )
    \end{pmatrix}.
\end{equation}
The permeability in the reference system $\tilde{\bm{\mu }}$ is given by the Voigt and Reuss limits with
\begin{equation}
    \tilde{\bm{\mu }} = \begin{pmatrix}
    \mu_{\mathrm{Voigt}} & 0\\
    0 & \mu _{\mathrm{Reuss}}
    \end{pmatrix},
\end{equation}
using (\ref{eq:modeling:multiscale:voigt2}) and (\ref{eq:modeling:multiscale:reuss2}).

\subsection{Numerical methods}

\subsubsection{Isogeometric Analysis}

In the context of the FEM, Isogeometric Analysis (IGA) is a technique which uses B-Splines and NURBS  to represent both geometry and solution of the partial differential equation (PDE). In the case of 2D magnetostatics, the magnetic vector potential $A_z$ is parametrized with $n$ ansatz functions $N_i(\mathbf{x})$ and coefficients $u_i$ as
\begin{equation}
    A_{z}(\mathbf{x}) = \sum\limits _{i=1}^{n} N_{i}(\mathbf{x}) u_{i}.
\end{equation}
To flexibly represent the permanent magnets, the same parametrization is  applied here to the remanence $\mathbf{B}_\mathrm{r}$, expressing  the $x$- and $y$-component as
\begin{align}
    \!B_{\mathrm{r}x}(\mathbf{x}) =\!\!\sum\limits _{i=1}^{n} N_{i}(\mathbf{x}) b_{i}^{( x)},~
    B_{\mathrm{r}y}(\mathbf{x}) =\!\!\sum\limits _{i=1}^{n} N_{i}(\mathbf{x}) b_{i}^{( y)}\!,\!
    \label{eq:Methodology:RemanenceDiscretization}
\end{align}
with coefficients $b_{i}^{( x)}$ and $b_{i}^{( y)}$, respectively.

Expressing both geometry and solution as NURBS and B-Splines makes IGA also a powerful tool for shape optimization, as (almost) arbitrary modifications of the geometry are possible by adjusting the coordinates of the control points without having to remesh the model \citep{ShapeOptimizationOfRotatingElectricMachines,IGAshapeOptimization}.

\subsubsection{Implementation of rotation}
In the context of IGA, harmonic mortaring has been proposed to efficiently simulate rotating electric machines, see, e.g.,  \citep{ShapeOptimizationOfRotatingElectricMachines, BONTINCK201840, OnTorqueComputationInElectricMachine}. This method is also suitable for the case of a magnetocaloric cooling device.
Instead of remeshing each rotor position, the boundaries of the rotating structures are coupled by harmonic (Fourier) basis functions. This allows for the use of non-conforming geometries at the coupling boundary which reduces the computational simulation effort. 

To this end, we define the domains $\OmegaRotor$ for the rotating part, $\OmegaStator$ for the stationary part and the boundary $\GammaAirGap$, which represents the interface in the air gap between the rotor and stator. In \cref{fig:Optimization:ShapeOpti:Initial}, $\GammaAirGap$ corresponds to the circle in the air gap between the AMRs and the rotor.
In the domains, the magnetostatic equation (\ref{eq:methodology:strong2DAnisotropic}) must be fulfilled, which simplifies to%
\begin{equation}
    \begin{cases}
    \nabla \cdotp \left( \mu _{\rt}^{-1} \nabla A_{z,\rt}\right) =\mu _{\rt}^{-1} \nabla \cdotp \mathbf{B}_{\mathrm{r}}^{\bot } & \mathrm{in} \ \Omega _{\rt}\\
    \nabla \cdotp \left( \mu _{\st}^{-1} \nabla A_{z,\st}\right) =\mu _{\st}^{-1} \nabla \cdotp \mathbf{B}_{\mathrm{r}}^{\bot } & \mathrm{in} \ \Omega _{\st},
    \end{cases} 
    \label{eq:NumericalMethods:Mortaring:StrongFormDomains}
\end{equation}
assuming isotropic material properties and no excitation currents. On the coupling boundary $\GammaAirGap$, the two conditions 
\begin{equation}
    \begin{cases}
    A_{z,\st}( \theta ) =A_{z,\rt}( \theta -\beta ) & \mathrm{on} \ \Gamma _{\ag}\\
    H_{\theta ,\st}( \theta ) = H_{\theta ,\rt}( \theta -\beta ) & \mathrm{on} \ \Gamma _{\ag},
    \end{cases}
    \label{eq:NumericalMethods:Mortaring:BoundaryConditions}
\end{equation}
must be fulfilled, where $\theta$ is the angle used for evaluation and $\beta$ the rotation angle. The field quantities are evaluated in their respective local coordinate system, see \citep{BONTINCK201840}. $H_\theta$ denotes the azimuthal component of the magnetic field strength and is given by
\begin{equation}
\begin{cases}
    H_{\theta ,\rt} = -\mu _{\rt}^{-1} \nabla A_{z,\rt} \cdotp \mathbf{n}_{\rt} \\
    H_{\theta ,\st} =\mu _{\st}^{-1} \nabla A_{z,\st} \cdotp \mathbf{n}_{\st} 
    \label{eq:NumericalMethods:Mortaring:AzimuthalComponents}
    \end{cases}
\end{equation}
with the outward pointing normal vectors $\mathbf{n}_\rt$ and $\mathbf{n}_\st$. Discretization of (\ref{eq:NumericalMethods:Mortaring:StrongFormDomains}) with the standard Ritz-Galerkin approach using B-Spline ansatz functions and the parametrization of (\ref{eq:NumericalMethods:Mortaring:BoundaryConditions}) and (\ref{eq:NumericalMethods:Mortaring:AzimuthalComponents}) using harmonic basis functions leads to the coupled system
\begin{equation}
    \begin{pmatrix}
    \mathbf{K}_{\rt} & 0 & -\mathbf{G}_{\rt}\\
    0 & \mathbf{K}_{\st} & \mathbf{G}_{\st}\\
    -\mathbf{G}_{\rt}^{\top} & \mathbf{G}_{\st}^{\top}& 0
    \end{pmatrix}\begin{pmatrix}
    \mathbf{u}_{\rt}\\
     \begin{array}{l}
    \mathbf{u}_{\st}\\
    \bm{\lambda }
    \end{array}
    \end{pmatrix} =\begin{pmatrix}
    \mathbf{b}_{\rt}\\
     \begin{array}{l}
    \mathbf{b}_{\st}\\
    0
    \end{array}
    \end{pmatrix},
    \label{eq:NumericalMethods:Mortaring:MatrixSystem}
\end{equation}
which is a saddle point problem as known from the literature \citep{Egger_2021}. Here, the stiffness matrix $\mathbf{K}$, the discrete solution vector for the magnetic vector potential $\mathbf{u}$, and the contribution of the permanent magnets $\mathbf{b}$ are separated on the domains $\OmegaRotor$ and $\OmegaStator$.
Further details on the coupling matrices $\mathbf{G}_{\rt}$, $\mathbf{G}_{\st}$ and the coupling coefficients $\bm{\lambda}$ can be found in \citep{OnTorqueComputationInElectricMachine, BONTINCK201840}.

Using this formulation, the torque can be directly computed using the Lagrange multipliers $\bm{\lambda}$, as explained in \citep{OnTorqueComputationInElectricMachine}.

\section{Structural optimization}
\label{chap:StructuralOptimization}
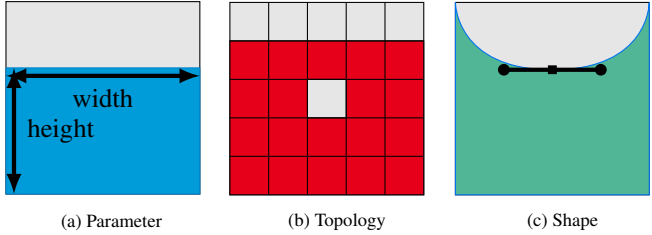
\begin{figure}[t]
    \centering
    \begin{subfigure}{.32\linewidth}
    \newcommand{\tw}{0.9\linewidth}
    \begin{tikzpicture}[y=1pt, x=1pt, yscale=-1, xscale=1, inner sep=0pt, outer sep=0pt]
					\draw[fill=black!10] (  0.0,  0.0) rectangle (\tw,\tw);
				    \draw[fill=TUDa-2a, draw=TUDa-2b] (  0.0, 25.0) rectangle (\tw,\tw);
					\draw[latex-latex, ultra thick]   (  3.0, 25.0) -- (  3.0,\tw) node [midway, right=0.5em] {height};
					\draw[latex-latex, ultra thick]   (  0.0, 28.0) -- (\tw, 28.0) node [midway, below=0.5em] {width};
				\end{tikzpicture}
    \caption{Parameter}
    \label{fig:Optimization:Parameter}
    \end{subfigure} \hfill
    \begin{subfigure}{.32\linewidth}
    \newcommand{\tw}{0.9\linewidth}
      \begin{tikzpicture}[y=1pt, x=1pt, yscale=-1, xscale=1, inner sep=0pt, outer sep=0pt]
			\draw[fill=black!10] (  0.0,  0.0) rectangle (\tw,\tw);
			\draw[fill=TUDa-9b ] (  0.0, 0.2*\tw) rectangle (\tw,\tw);
			\draw[fill=black!10] ( 0.4*\tw, 0.4*\tw) rectangle ( 0.6*\tw, 0.6*\tw);
			\foreach \i in {1,...,4}
			{
				\draw[] (\i*0.2*\tw,0) -- (\i*0.2*\tw,\tw);
				\draw[] (0,\i*0.2*\tw) -- (\tw,\i*0.2*\tw);
			}
		\end{tikzpicture}
		\caption{Topology}
		\label{fig:Optimization:Topology}
    \end{subfigure} \hfill
    \begin{subfigure}{.32\linewidth}
    \newcommand{\tw}{0.9\linewidth}
      \begin{tikzpicture}[y=1pt, x=1pt, yscale=-1, xscale=1, inner sep=0pt, outer sep=0pt]
			\draw[fill=black!10] (  0.0,  0.0) rectangle (\tw,\tw);
			\draw[TUDa-3d,fill=TUDa-3a] (  0.0,\tw) -- (  0.0,  0.0) .. controls (  0.0,0.0000) and (  0.0,0.35*\tw) .. ( 0.5*\tw,0.35*\tw) .. controls (\tw,0.35*\tw) and (\tw,0.0000) .. (\tw,  0.0) -- (\tw,\tw) -- cycle;
			\draw[fill=black] ( 0.25*\tw, 0.35*\tw) circle (2);
			\draw[fill=black] ( 0.75*\tw, 0.35*\tw) circle (2);
			\draw[fill=black] (0.48*\tw,0.33*\tw) rectangle (0.52*\tw,0.37*\tw);
			\draw[ultra thick,black] ( 0.25*\tw, 0.35*\tw) -- ( 0.75*\tw, 0.35*\tw);
		\end{tikzpicture}
		\caption{Shape}
		\label{fig:Optimization:Shape}
    \end{subfigure} \hfill
    \caption{Types of structural optimization.}
    \label{fig:Optimization:Types}
\end{figure}

Structural optimization is typically classified into the three types \textit{parameter}, \textit{topology} and \textit{shape} optimization \citep{harzheim2008strukturoptimierung}. For a given geometry, parameter optimization is used to determine optimal (geometric) parameters, such as width, height or radius (\cref{fig:Optimization:Parameter}). In the case of topology optimization, the aim is to find novel designs, e.g., by assigning local (pseudo) density values for solid and void regions. This allows for adding or removing holes in the geometry (\cref{fig:Optimization:Topology}). The result serves as draft for a specific structural setting.  Given the basic topology, shape optimization can be used to adjust the boundary of the geometry. In the case of IGA, this can be achieved by moving control points that define the boundary (\cref{fig:Optimization:Shape}).

In this work, topology optimization is first used to determine the best topological layout. Shape optimization is then applied to further homogenize the magnetic field. 

\subsection{Topology optimization}
Research concerning the topology of magnet assemblies to date usually assumes a fixed remanence direction given by Halbach cylinders \citep{PermanentMagnetDesignForMagneticHeatPumps, TopologyOptimizedPermanentMagnetSystems}. 
This work proposes a more general design approach, where the local remanence may change its magnetization direction, allowing for novel designs to be created.

\subsubsection{Methodology}

Larger changes in the magnetic field induce higher temperature changes in the MCM. Hence, it is necessary to have a high magnetic flux density in the air gap during magnetization and a near-zero magnetic flux density during demagnetization. Consequently, the objective function is chosen to maximize the difference of the average magnetic flux density during magnetization $\langle \Bhigh \rangle$ and the average magnetic flux density during demagnetization $\langle \Blow \rangle$.  
Here, a symmetric operation of the AMRs for a two-pole refrigerator is assumed, i.e., the (de)magnetization phases are equally sized to simplify operation and assembly of the system.  

\begin{figure}[t]
    \centering
    \begin{subfigure}{.48\linewidth}
      \centering
        \includegraphics[trim= 40 27 23 8, clip, width=\linewidth, height=\linewidth]{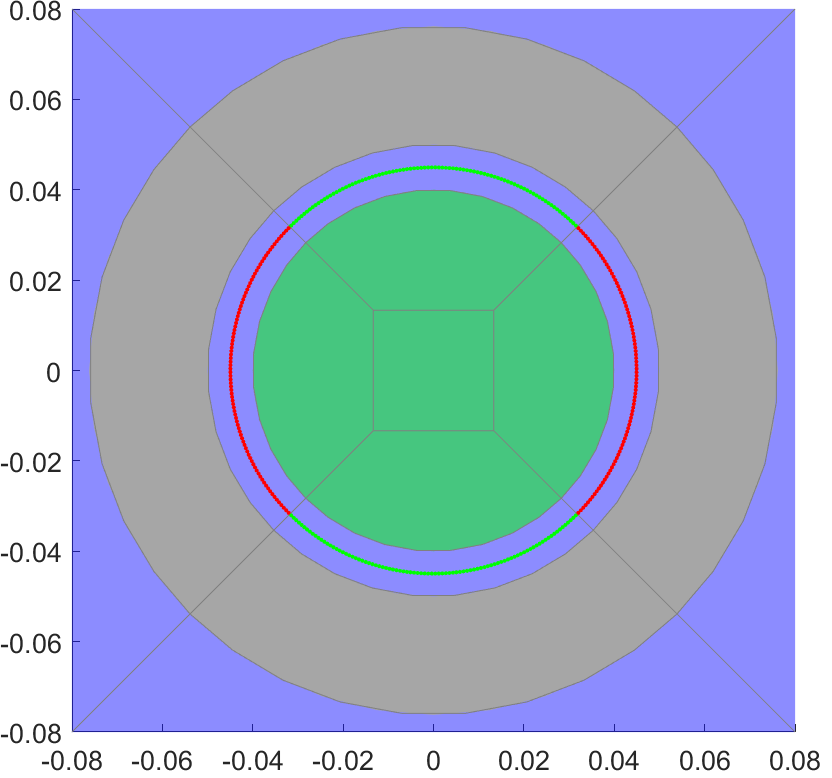}
        \caption{Inner rotary case}
      \label{fig:optimization:designArea1}
    \end{subfigure} \hfill
    \begin{subfigure}{.48\linewidth}
      \centering
        \includegraphics[trim= 45 27 20 10, clip, width=\linewidth, height=\linewidth]{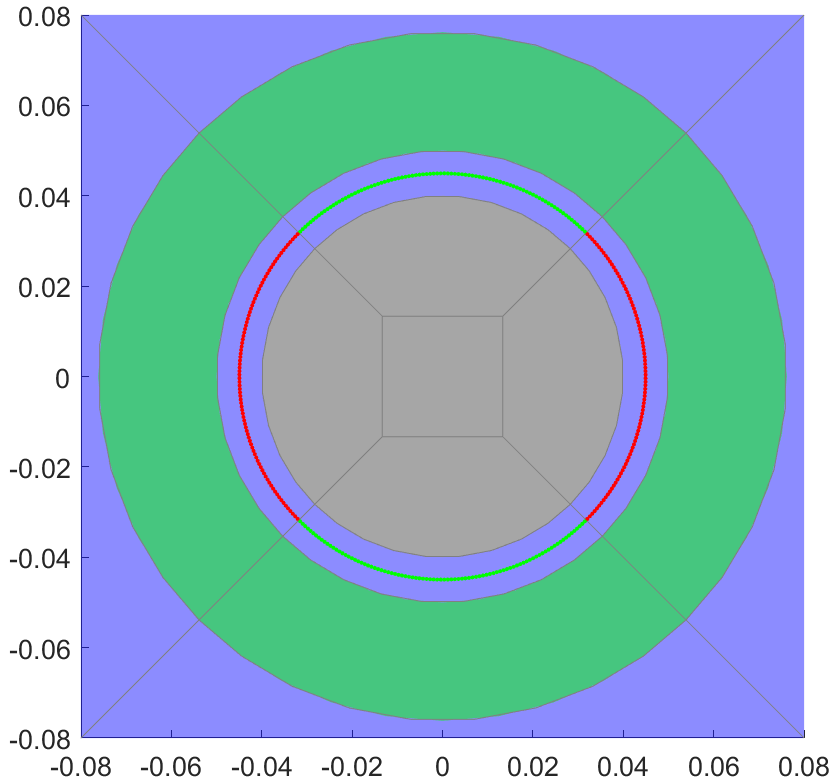}
      \caption{Outer rotary case}
      \label{fig:optimization:designArea2}
    \end{subfigure}
    \caption{Design areas for the inner and outer rotary case. The green patches indicate the design area for the topology optimization. $\langle \Bhigh \rangle$ and $\langle \Blow \rangle$ are obtained by evaluation of the indicated red and green points.}
    \label{fig:optimization:designAreas}
\end{figure}

Moreover, only rotating permanent magnets are investigated, as rotating AMRs have shown disadvantages regarding the hydraulic system and sealing friction~\citep{ExperimentalResultsForANovelRotaryAMR}. The design areas for the topology optimization are shown as green patches in the Figs. \ref{fig:optimization:designArea1} (inner rotary) and \ref{fig:optimization:designArea2} (outer rotary). Gray and blue patches represent iron and air, respectively. 
The combination of both design areas leads to the co-rotational case. The magnetic flux density in the air gap is evaluated at equidistant points indicated in green ($\Bhigh$) and red ($\Blow$). For the inner radius, air gap thickness and outer radius, $r_\mathrm{inner}=\SI{40}{\mm}$, $d_\mathrm{ag}=\SI{10}{\mm}$ and $r_\mathrm{outer}=\SI{75}{\mm}$ are used, respectively. However, these values can be scaled without influencing the distribution of the  magnetic flux density due to the scale independence of magnetic fields in the absence of currents or charges \citep{coey_2010}.

The remanence $\mathbf{B}_\mathrm{r}$ in the design area is discretized according to (\ref{eq:Methodology:RemanenceDiscretization}). The coefficients $b_i^{(x)}$ and $b_i^{(y)}$ are the design variables and will be modified during the optimization. Depending on the absolute value of $\mathbf{B}_\mathrm{r}$, the material behavior resembles either a permanent magnet or soft iron. 

Two additional constraints need to be imposed.
First, the remanence is limited to a physically feasible value $B_\mathrm{r,mag}$. Second, the remanent area, which correlates with the volume of the permanent magnet, is restricted by $A_{\mathrm{max}}$ to allow for an effective comparison between different designs. The resulting PDE-constrained optimization problem reads
\begin{align}
    \max_{\mathbf{B}_\mathrm{r}}  \quad & f_{\opt} = \langle B_{\mathrm{high}} \rangle  - \langle B_{\mathrm{low}} \rangle 
    \label{eq:Optimization:TopologyProblemMax}
    \\
    \textrm{s.t.} \quad & \norm{\mathbf{B}_\mathrm{r}}  \ \leq \ B_{\mathrm{r,mag}}, \\
      & \int _{\Omega } \norm{\mathbf{B}_\mathrm{r}} \Integ \ \leq \ B_{\mathrm{r,mag}} A_{\mathrm{max}},    \\
      \mathrm{and} \quad &  
      (\ref{eq:methodology:strong2DAnisotropic}) \mathrm{~discretized~with~IGA.}
    \label{eq:Optimization:TopologyProblemDefinition}
\end{align}

\begin{figure}[t]
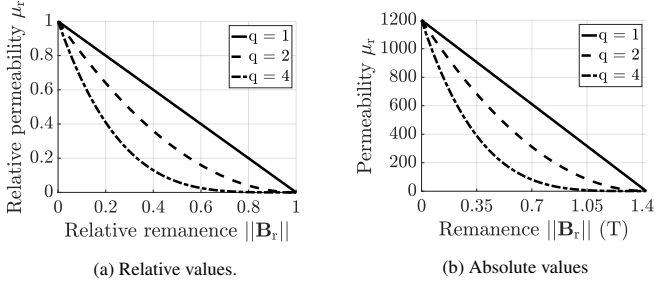

    \centering
    \begin{subfigure}{.48\linewidth}
      \centering
        \includesvg[width=\linewidth]{SIMPrel}
        \caption{Relative values.}
      \label{fig:optimization:simprel}
    \end{subfigure} \hfill
    \begin{subfigure}{.48\linewidth}
      \centering
        \includesvg[width=\linewidth]{SIMP}
      \caption{Absolute values}
      \label{fig:optimization:simpabs}
    \end{subfigure}
    \caption{SIMP relationship to couple remanence and permeability to allow for the continuous transition from soft iron to permanent magnet regions. Intermediate values are penalized with increasing $q$.}
    \label{fig:optimization:simp}
\end{figure}

In mechanical topology optimization, the elastic modulus is coupled by a Solid Isotropic Material with Penalization (SIMP) relationship to describe the transition from solid to void \citep{sigmund2013}. In analogy, an adapted SIMP approach is used in this work to couple remanence and permeability of the material to distinguish between magnet and soft iron. This is possible, as soft iron has a high permeability and no remanence, whereas permanent magnets have low permeability and high remanence. The relationship is modeled as a decreasing function
\begin{equation}
    \mu _{\mathrm{rel}} =( 1-\| \mathbf{B}_\mathrm{r} \| _{\mathrm{rel}})^{q}, \label{eq:TopOpt:MuRelSimp}
\end{equation} 
where $q$ is the penalization factor to avoid intermediate values of the material density in the topology optimization. 
The relative values for the remanence $\| \mathbf{B}_\mathrm{r} \| _{\mathrm{rel}} \in [0,1]$  and the relative permeability $\mu _{\mathrm{rel}}\in [0,1]$ are obtained by a linear mapping between the minimum and maximum values, which are given by the magnet and soft iron. Relative permeabilities of $\mu_{\mathrm{r},\mathrm{iron}}=\SI{1200}{}$ and $\mu_{\mathrm{r},\mathrm{pm}}=\SI{1.05}{}$ are used for soft iron and permanent magnet, respectively. The remanence limits are set to zero for iron and $|B_\mathrm{r,mag}| = \SI{1.425}{\tesla}$ for the magnet, which is achievable, e.g., with Nd-Fe-B magnets. This dependence is demonstrated for different penalization factors in \cref{fig:optimization:simp}.  Here, a value of $q=4$ has shown to be appropriate. 

\subsubsection{Results} %
The implementations are made in MATLAB\textsuperscript{\textregistered} using the open source packages \textit{NURBS} \citep{nurbsToolbox}, \textit{Tensor Toolbox} \citep{tensorToolbox} and \textit{GeoPDEs} \citep{geopdes, geopdes3.0} as well as the MATLAB\textsuperscript{\textregistered} \textit{Optimization} toolbox \citep{matlabOptimization}. 

All optimization steps are performed with sensitivity analysis and analytical derivations of the gradient information.
The adjoint method \citep{harzheim2008strukturoptimierung} is used to derive the gradient of the objective function and constraints with respect to the remanence coefficients $b{_{i}^{(x)}}$ and $b{_{i}^{(y)}}$. 
The optimization is then realized with the \textit{interior-point} method from the MATLAB\textsuperscript{\textregistered} solver \textit{fmincon}. For comparability, 100 iterations are performed in each optimization. This achieves efficient computation and convergence within 30 minutes for each optimization on a PC with an Intel Dual Core i7-4600U processor and 8 GB of RAM.

\begin{figure*}
  \centering
  \begin{minipage}{\linewidth}
    \begin{minipage}{0.245\linewidth}
        \centering
        \subcaptionbox{Iteration 0.\label{fig:TopOpt:Results:InnerRotating:TopBottom:Iter0}}
          {\includegraphics[trim= 235 119 265 102, clip, width=\linewidth, height=\linewidth]{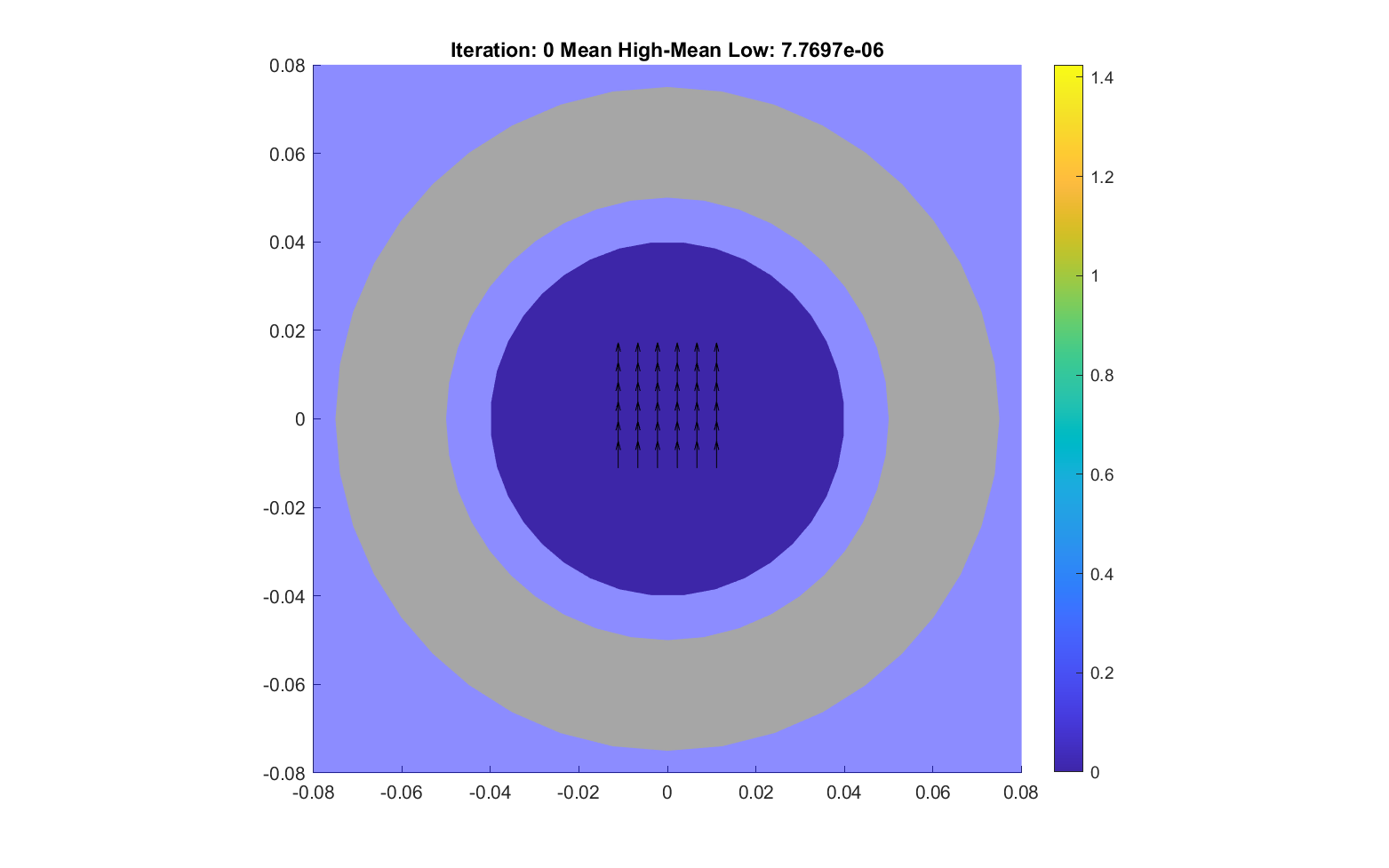}}
    \end{minipage}
    \begin{minipage}{0.245\linewidth}
        \centering
        \subcaptionbox{Iteration 10.\label{fig:TopOpt:Results:InnerRotating:TopBottom:Iter10}}
          {\includegraphics[trim= 235 119 265 102, clip, width=\linewidth, height=\linewidth]{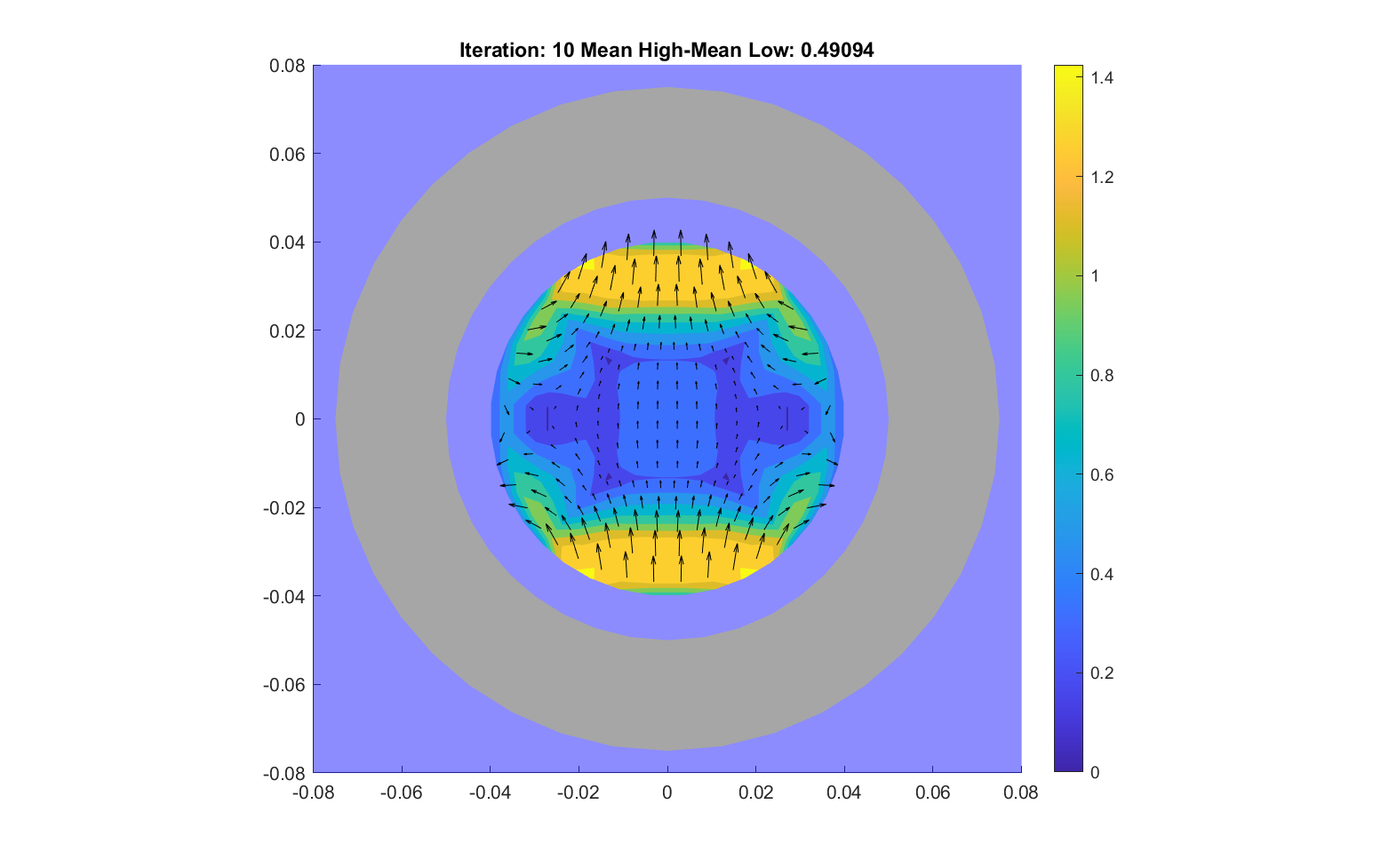}}
    \end{minipage}
    \begin{minipage}{0.245\linewidth}
        \centering
        \subcaptionbox{Iteration 50.\label{fig:TopOpt:Results:InnerRotating:TopBottom:Iter50}}
          {\includegraphics[trim= 235 119 265 102, clip, width=\linewidth, height=\linewidth]{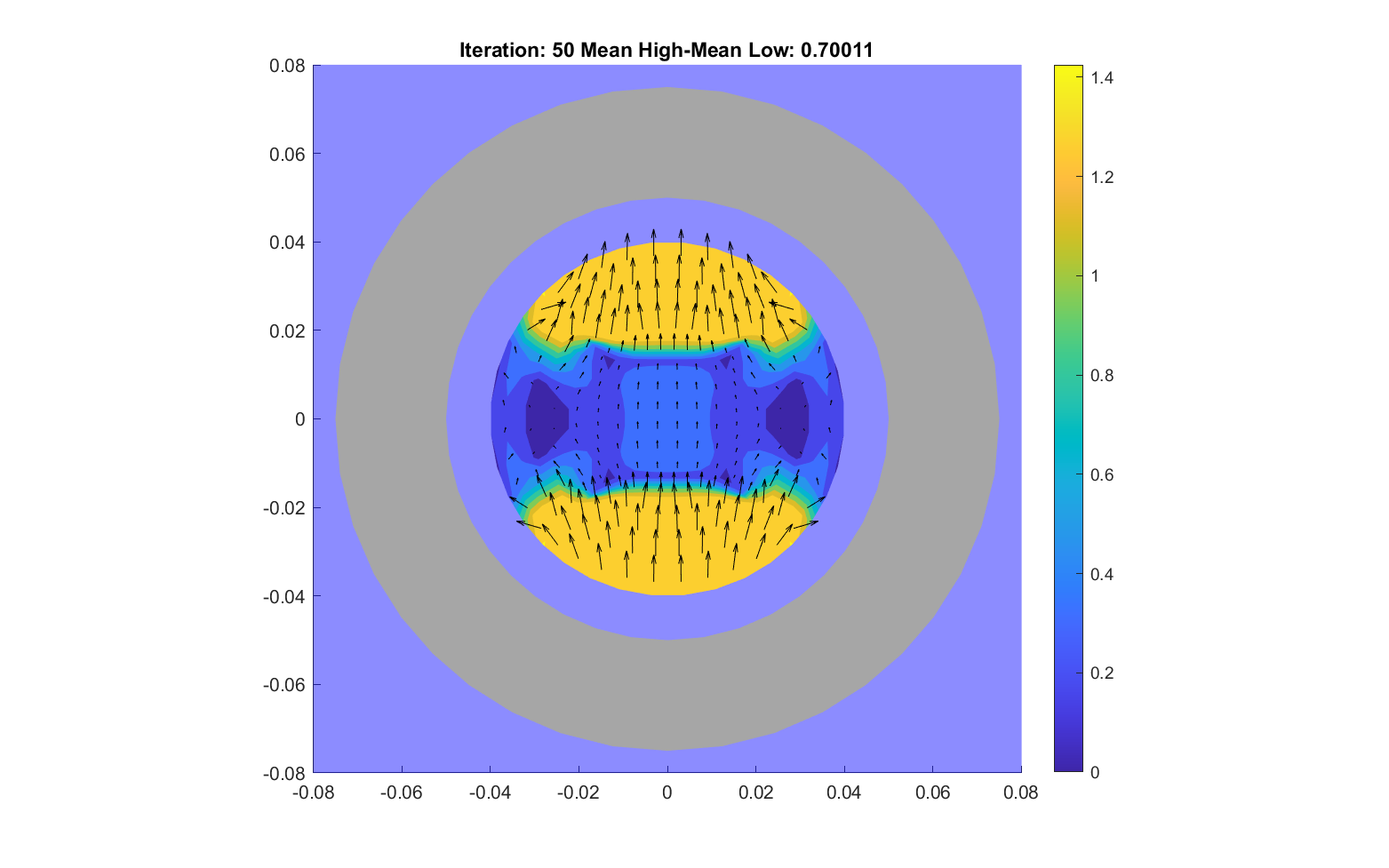}}
    \end{minipage}
    \begin{minipage}{0.245\linewidth}
        \centering
        \subcaptionbox{Iteration 100.\label{fig:TopOpt:Results:InnerRotating:TopBottom:Iter100}}
          {\includegraphics[trim= 235 119 265 102, clip, width=\linewidth, height=\linewidth]{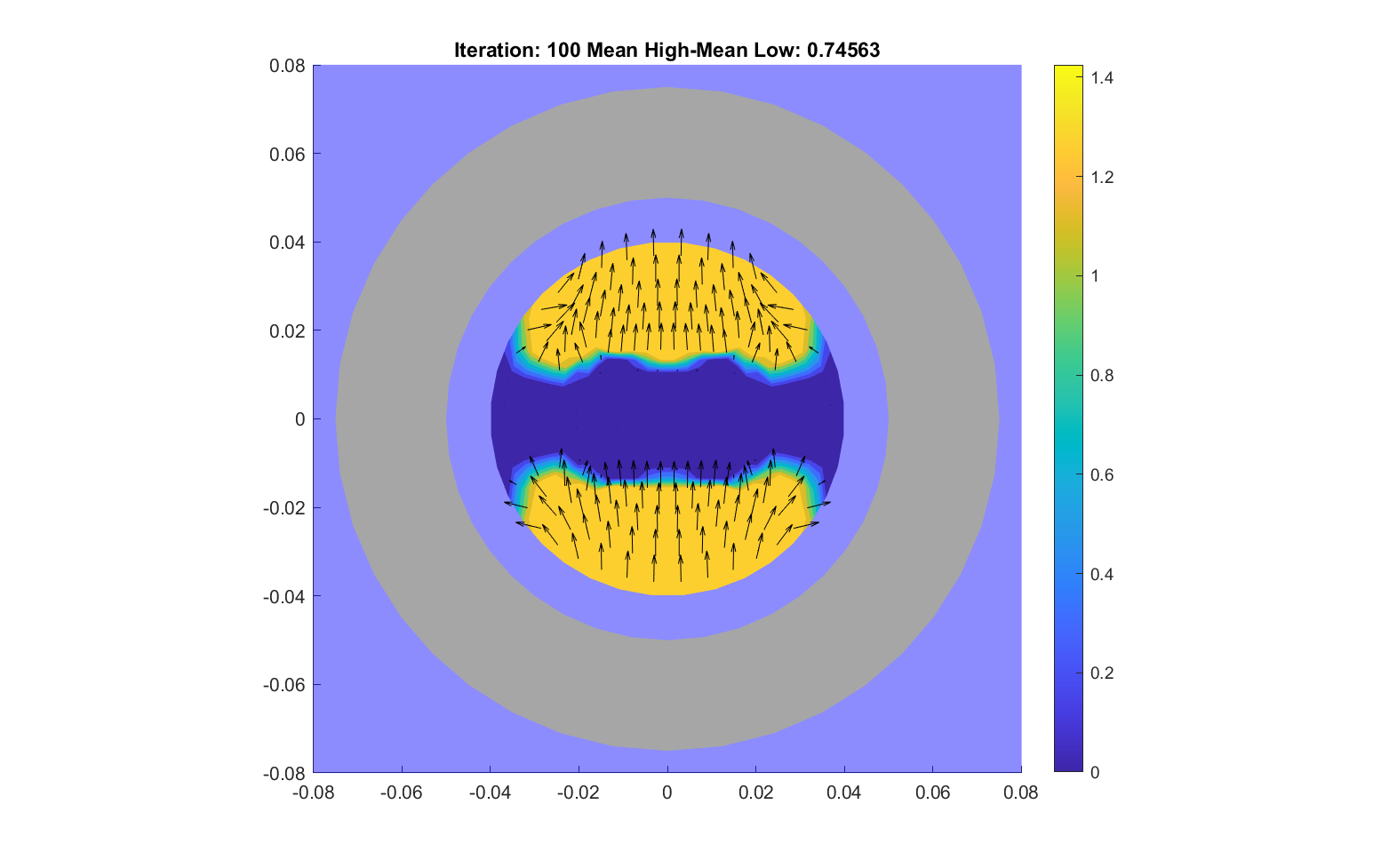}}
    \end{minipage}

    \caption{Optimization steps for the top and bottom inner rotary design. Yellow areas indicate permanent magnet material, deep blue and gray areas represent soft iron. A low $B_\mathrm{r}$ of \SI{0.02}{\tesla} in the center patch prescribes the magnetization direction in iteration 0. The areas directly aligning with the air gap have the
highest contribution to the objective function and are increased in remanence, leading to a separated top and bottom design. }
    \label{fig:TopOpt:Results:InnerRotating:TopBottom}
    
    \bigskip
      \begin{minipage}{0.245\linewidth}
        \centering
        \subcaptionbox{Iteration 0.\label{fig:TopOpt:Results:InnerRotating:Centered:Iter0}}
          {\includegraphics[trim= 239 119 263 100, clip, width=\linewidth, height=\linewidth]{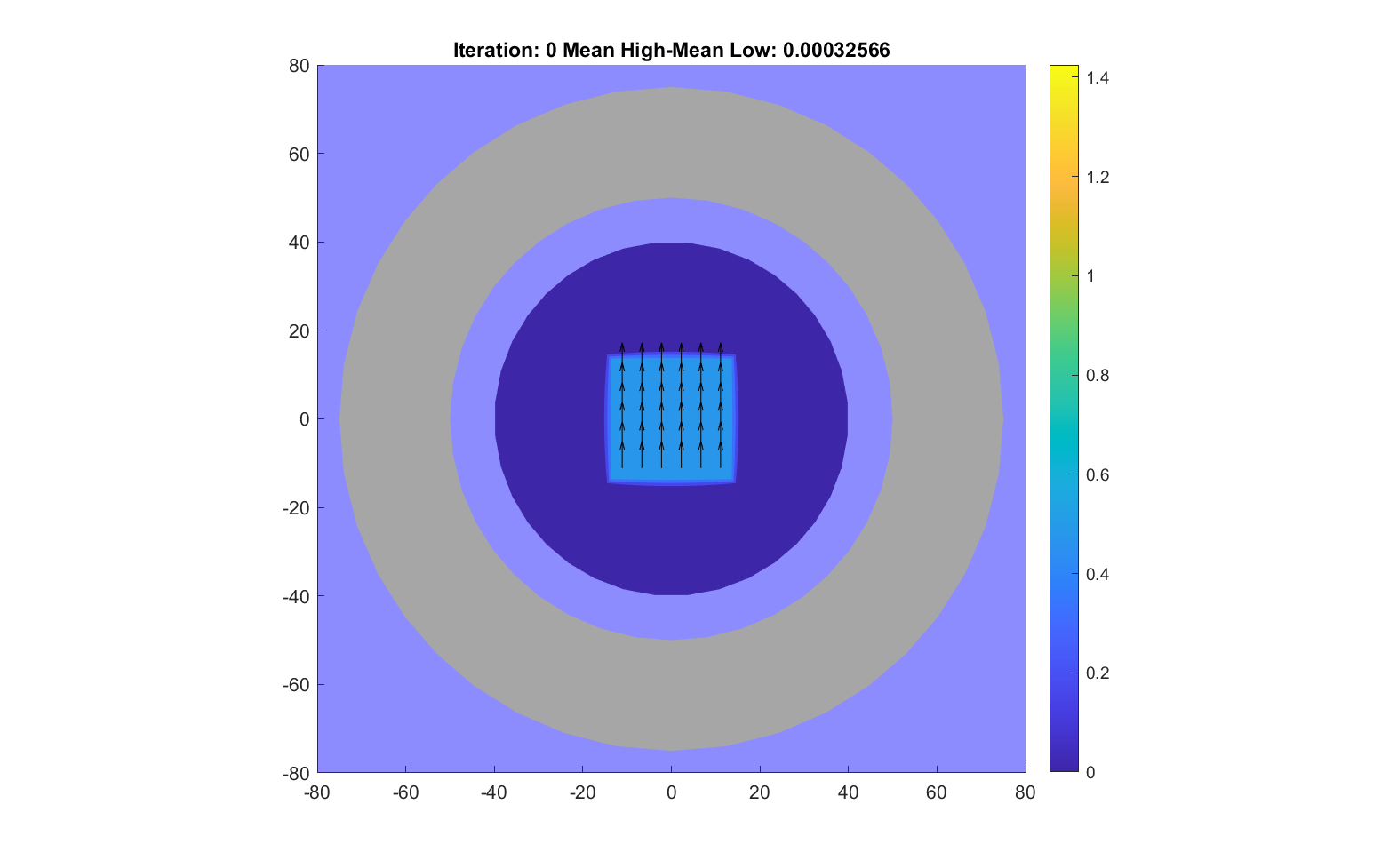}}
    \end{minipage}
    \begin{minipage}{0.245\linewidth}
        \centering
        \subcaptionbox{Iteration 10.\label{fig:TopOpt:Results:InnerRotating:Centered:Iter10}}
          {\includegraphics[trim= 239 119 263 100, clip, width=\linewidth, height=\linewidth]{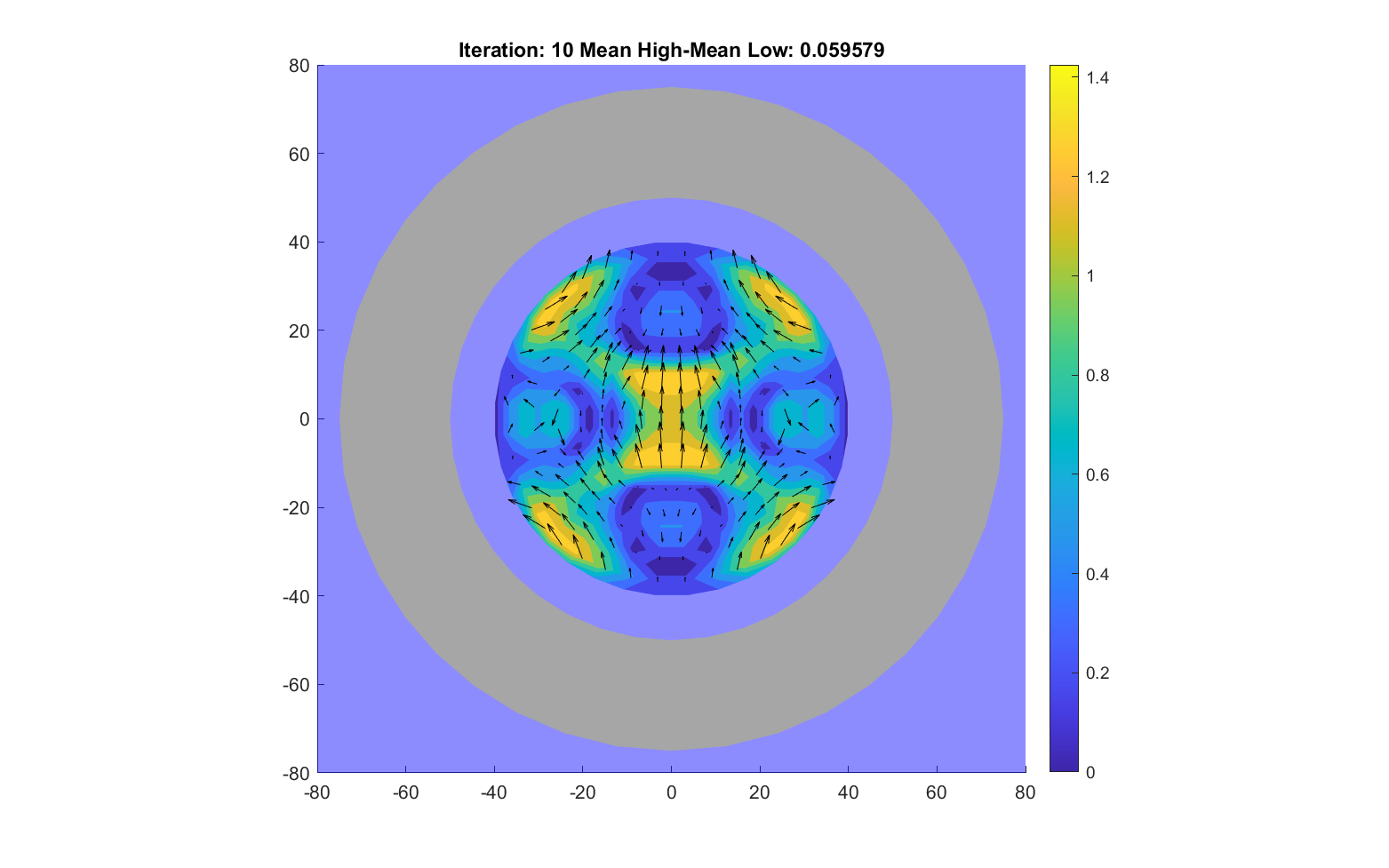}}
    \end{minipage}
    \begin{minipage}{0.245\linewidth}
        \centering
        \subcaptionbox{Iteration 50.\label{fig:TopOpt:Results:InnerRotating:Centered:Iter50}}
          {\includegraphics[trim= 239 119 263 100, clip, width=\linewidth, height=\linewidth]{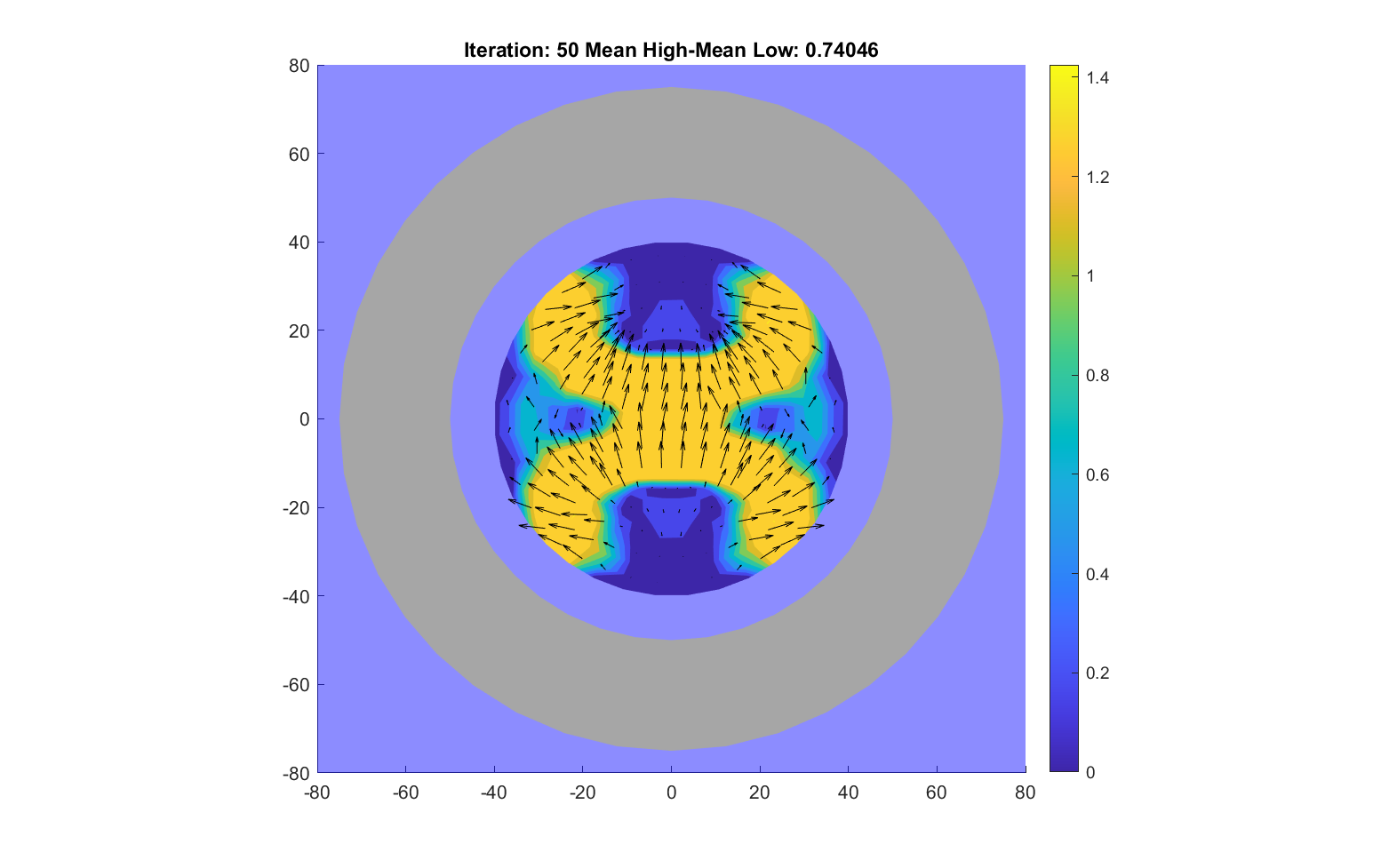}}
    \end{minipage}
    \begin{minipage}{0.245\linewidth}
        \centering
        \subcaptionbox{Iteration 100.\label{fig:TopOpt:Results:InnerRotating:Centered:Iter100}}
          {\includegraphics[trim= 239 119 263 100, clip, width=\linewidth, height=\linewidth]{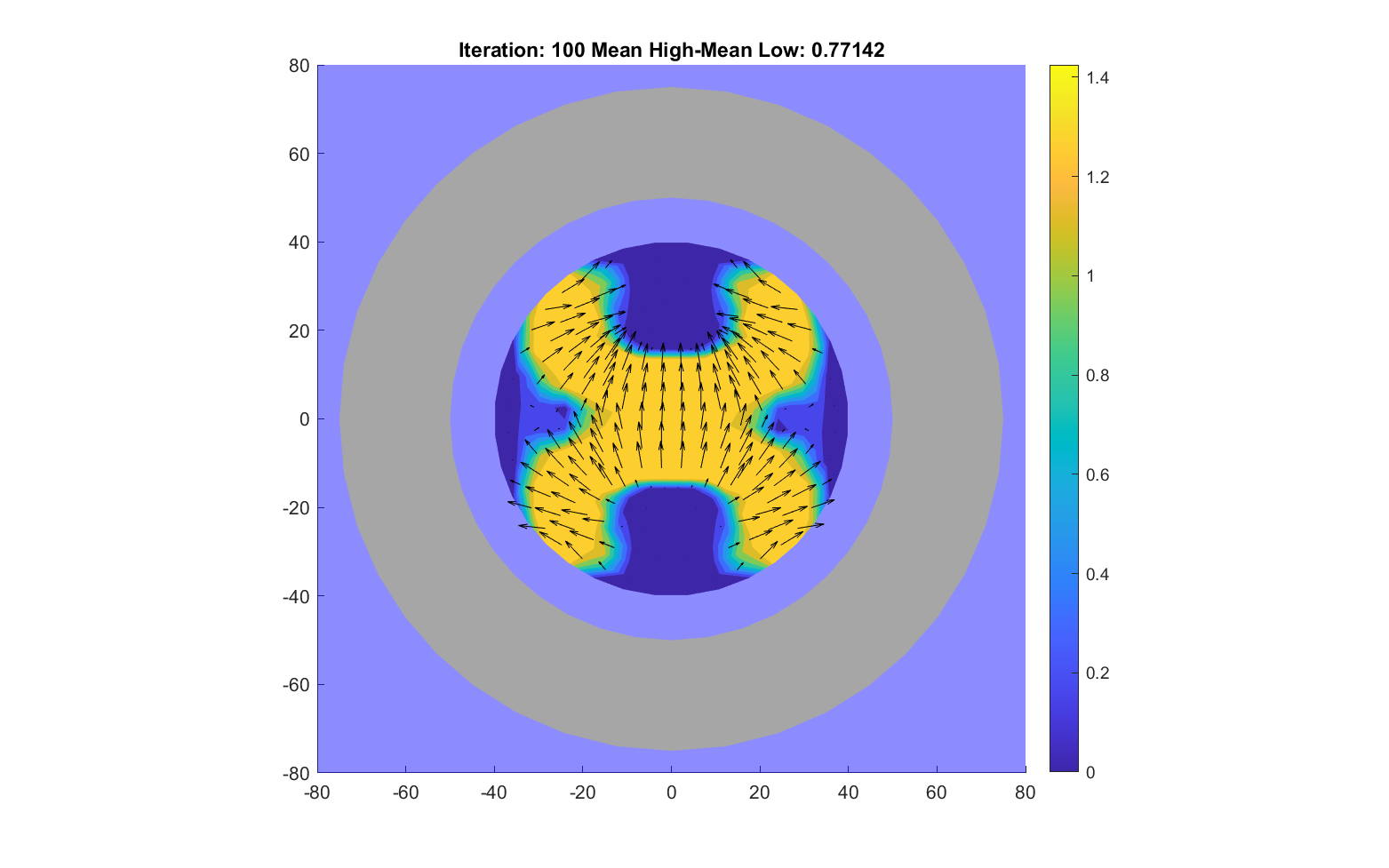}}
    \end{minipage}

    \caption{Optimization steps for the centered inner rotary design. Yellow areas indicate permanent magnet material, deep blue and gray areas represent soft iron. For a higher $B_\mathrm{r}$ of \SI{0.85}{\tesla} as initial condition for the center patch, the solver converges to a centered design. The already remanent area is used by adding further remanent material at the patch corners. The magnetic flux is then concentrated in the areas, which take the properties of soft iron.}
    \label{fig:TopOpt:Results:InnerRotating:Centered}

  \bigskip
   \begin{minipage}{0.245\linewidth}
        \centering
        \subcaptionbox{Inner rotary top and bottom.\label{fig:TopOpt:Results:InnerRotating:TopBotFinal}}
          {\includegraphics[trim= 176 59 206 40, clip, width=\linewidth, height=\linewidth]        {L1Remanence100.png}}
    \end{minipage}
    \begin{minipage}{0.245\linewidth}
        \centering
        \subcaptionbox{Inner rotary centered.\label{fig:TopOpt:Results:InnerRotating:CenteredFinal}}
          {\includegraphics[trim= 179 59 203 40, clip, width=\linewidth, height=\linewidth] {L2Remanence100.png}}
    \end{minipage}
    \begin{minipage}{0.245\linewidth}
        \centering
        \subcaptionbox{Outer rotary.\label{fig:TopOpt:Results:OuterRotating:Final}}
          {\includegraphics[trim= 176 59 206 40, clip, width=\linewidth, height=\linewidth]{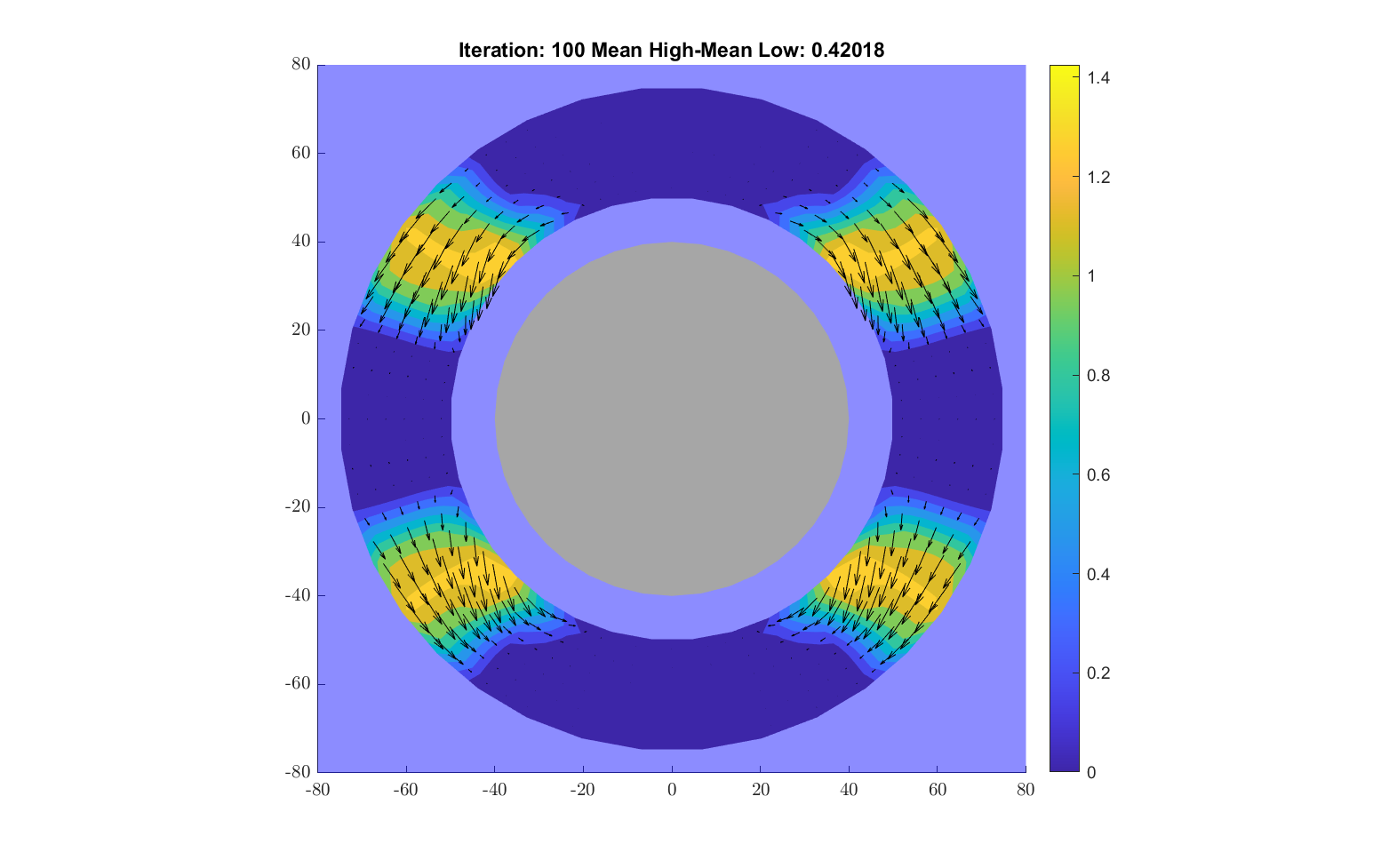}}
    \end{minipage}
    \begin{minipage}{0.245\linewidth}
        \centering
        \subcaptionbox{Co-rotational.\label{fig:TopOpt:Results:CoRotatingFinal}}
          {\includegraphics[trim= 176 59 206 40, clip, width=\linewidth, height=\linewidth]{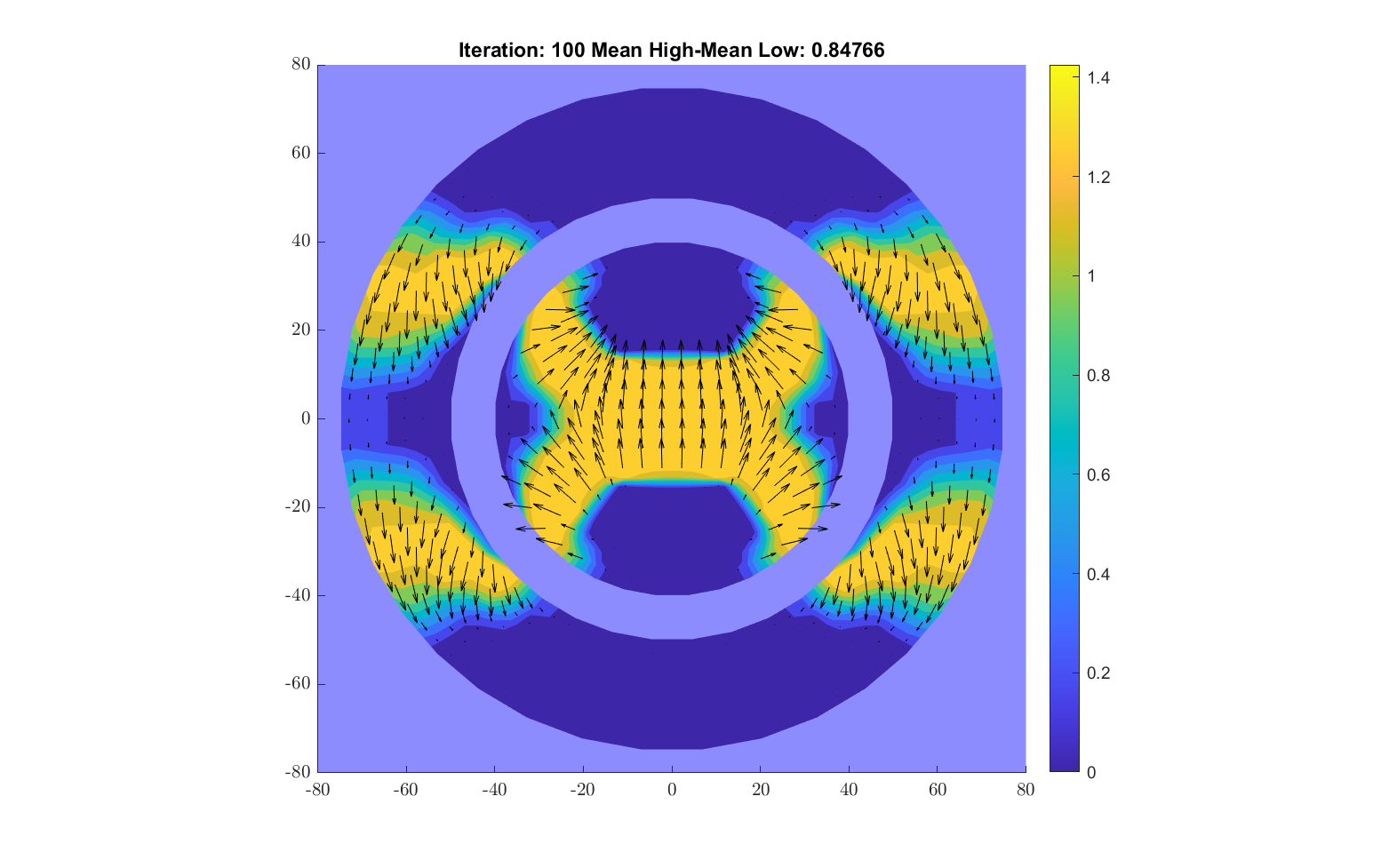}}
    \end{minipage}

    \caption{Topology optimization results for the inner rotary, outer rotary and co-rotational case. Yellow areas indicate permanent magnet material, deep blue and gray areas represent soft iron. For the inner rotary case, two solutions are found, depending on the initial conditions (\cref{fig:TopOpt:Results:InnerRotating:TopBotFinal} and \ref{fig:TopOpt:Results:InnerRotating:CenteredFinal}). In the outer rotary case (\cref{fig:TopOpt:Results:OuterRotating:Final}), permanent magnet material is positioned at the transition zones from magnetization to demagnetization to generate a sharper distinction for high and low magnetic field areas. The co-rotational case (\cref{fig:TopOpt:Results:CoRotatingFinal}) combines the inner and outer rotary case (also depending on initial conditions). }
    \label{fig:TopOpt:Results:Final}
  \end{minipage}
\end{figure*}

Due to the nature of gradient based algorithms, the initial conditions influence the  convergence process and hence the final solution. This is demonstrated for the inner rotary design in the Figures \ref{fig:TopOpt:Results:InnerRotating:TopBottom} and \ref{fig:TopOpt:Results:InnerRotating:Centered}. The dark blue areas indicate non-remanent regions, which behave like soft iron. Yellow areas are assigned the properties of a permanent magnet, where the remanence direction is given by the black arrows. Intermediate values are plotted in the color range in between.

To magnetize the air gap, the remanence can be either directed from bottom to top or vice versa. Due to this symmetry, the gradient is zero if all-zero initial conditions are provided, as this corresponds to a saddle point of the objective function. In \cref{fig:TopOpt:Results:InnerRotating:TopBottom:Iter0}, a small value is assigned to the remanence of the center patch. As the patches aligning with the air gap have the highest contribution to the objective function, the remanence is increased exactly there during the iterations (\cref{fig:TopOpt:Results:InnerRotating:TopBottom:Iter10}). During the optimization process, intermediate values are removed and more material becomes remanent (\cref{fig:TopOpt:Results:InnerRotating:TopBottom:Iter50}), until the specified maximum magnetic area is reached (\cref{fig:TopOpt:Results:InnerRotating:TopBottom:Iter100}), resulting in the \textit{inner rotary top and bottom} design, similar to the one in \citep{Dimitrov2006aa}. For its realization, sintered magnets are needed, where the magnetization is uniformly approximated.

When using slightly different initial conditions with a remanent value of $| \mathbf{B}_\mathrm{r}|=\SI{0.85}{\tesla}$ (Fig \ref{fig:TopOpt:Results:InnerRotating:Centered:Iter0}), the solver converges to a second possible design, the \textit{inner rotary centered} design. 
Increasing the remanence at the areas of the air gap in this case would also decrease the permeability there and therefore impede the already existing magnetic flux. Instead, creating a remanent arc (\cref{fig:TopOpt:Results:InnerRotating:Centered:Iter10})  increases the objective function faster. After several iterations (Fig \ref{fig:TopOpt:Results:InnerRotating:Centered:Iter50}), the centered design is obtained (\cref{fig:TopOpt:Results:InnerRotating:Centered:Iter100}). Building it  requires multiple magnet segments and the omission of a continuous shaft, making the assembly complex. At the same time, replacing the shaft is crucial to remove the parasitic backflow of magnetic flux through the shaft, seen in designs from the literature \citep{Kanluang_Simulations}.

\begin{figure}
    \centering
    \includesvg[width=\linewidth]{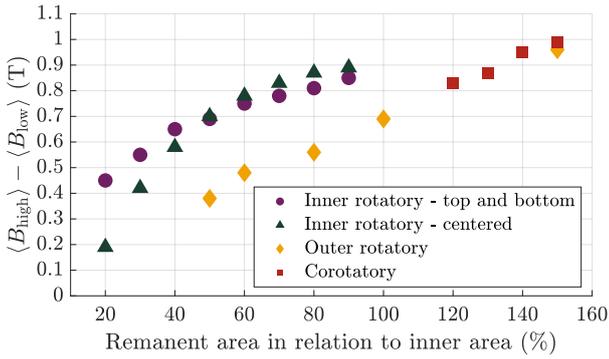}
    \caption{
    Comparison of the inner, outer and co-rotational design variants. The objective value (\ref{eq:Optimization:TopologyProblemMax}) is compared for the four design variants, given in \cref{fig:TopOpt:Results:Final}. The inner design area from \cref{fig:optimization:designArea1} is used as reference for the comparison.
    }
    \label{fig:Optimization:TopologyComparison}
\end{figure}

Carrying out this optimization process also for the outer- and co-rotational cases yields four different design possibilities, given in \cref{fig:TopOpt:Results:Final}. 
In the \textit{outer rotary} design (\cref{fig:TopOpt:Results:OuterRotating:Final}), permanent magnet material is positioned at the transition zones from magnetization to demagnetization to generate a sharper distinction for the high and low magnetic field areas. The \textit{co-rotational} design (\cref{fig:TopOpt:Results:CoRotatingFinal}) combines the inner and outer rotary cases.

Repeating this optimization process for varying permissible permanent magnet volume or -- in 2D -- cross-section area $A_{\mathrm{max}}$ with the four design variants yields the results presented in \cref{fig:Optimization:TopologyComparison}. 
Note that the $x$-axis shows the relative area $A_{\%}=A_{\mathrm{max}}/A_{\mathrm{inner}}$ of the permanent magnets in relation to the inner area $A_{\mathrm{inner}}$ from \cref{fig:optimization:designArea1} to enable a comparison of the designs. 
Summarizing the most important findings from \cref{fig:Optimization:TopologyComparison}, we see:
\begin{itemize}
    \item For low magnetic areas, the \textit{inner rotary top and bottom} design performs best. If
the magnet area is severely restricted, putting the magnets exactly where needed yields the best results.
    \item For magnetic areas of $A_\%\geq\SI{50}{\percent}$, the \textit{inner rotary centered} design outperforms the \textit{top and bottom} one. With higher proportion of permanent magnets, the concentration of the magnetic flux proves advantageous. 
    \item The \textit{outer rotary} design yields the lowest objective values for all cases, as explained by the unfavorable scaling in the radial direction. 
    \item The overall highest objective values are reached with a \textit{co-rotational} design, as to be expected. It is important to note that the co-rotational case is not a superposition of the inner and outer rotary cases due to the low permeability of permanent magnets, which impedes the magnetic flux. Such high magnetic fields therefore come at the cost of a significantly higher permanent magnet volume.
\end{itemize}

When applying the optimization method to an adapted problem with a different number of poles or a large outer diameter, further interesting results are found:
\begin{itemize}
    \item Some articles feature a design with four magnetic poles \citep{YOU2016231,optimizedrefrigeration}. This has been tested, but the objective values are found to be \SIrange{15}{30}{\percent} smaller compared to a two pole design. It is therefore recommended to prefer two pole systems operated at higher frequency.
    \item For large outer diameters, an outer rotary magnet design with magnets adjacent to the air gap similar to \citep{OptimizingMagnetostaticAssemblies,Insigna2016} is found (top right yellow diamond in \cref{fig:Optimization:TopologyComparison}). However, such an assembly is discouraged for commercial applications due to high space and permanent magnet mass requirements.
\end{itemize}

\subsection{Design implementation}

Assuming that geometries similar to the ones obtained from the topology optimization yield a comparable magnetic flux density, the \textit{inner rotating centered} design (\cref{fig:TopOpt:Results:InnerRotating:CenteredFinal}) is chosen as most suitable initial topology for further optimization.

The \textit{inner rotary top and bottom} design is more efficient at low magnet mass, but with current MCM technology, higher magnetic fields are needed. 
The \textit{co-rotational} design yields higher objective values, but this is at the expense of a drastically larger magnet mass and a higher complexity of the assembly. This design is only recommended for applications where a higher power density is necessary.
Overall, the \textit{inner rotary centered} design is a good compromise between a high magnetic flux density and comparatively low permanent magnet mass. If needed, this design can be extended to a co-rotational one.

\begin{figure*}
    \centering
    \begin{subfigure}[t]{.48\linewidth}
      \centering
        \includegraphics[trim= 130 30 120 25, clip, width=\linewidth]{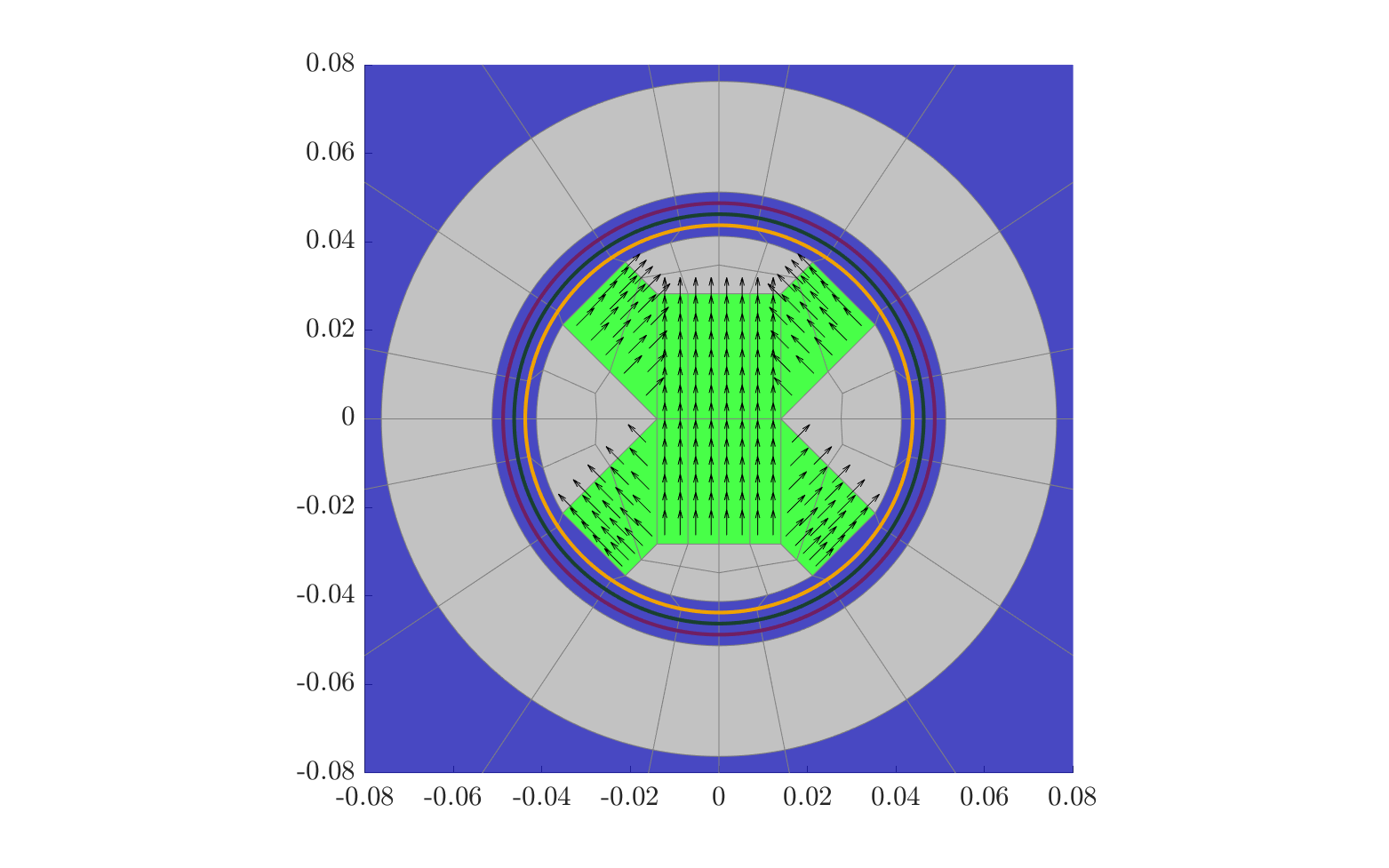}
        \caption{Geometry based on the topology optimization for the inner rotary centered design from \cref{fig:TopOpt:Results:InnerRotating:CenteredFinal}.  The lines in the air gap are used to evaluate the magnetic flux density.}
      \label{fig:optimization:resultsGeometry}
    \end{subfigure} \hfill
    \begin{subfigure}[t]{.48\linewidth}
      \centering
        \includegraphics[trim= 130 30 120 25, clip, width=\linewidth]{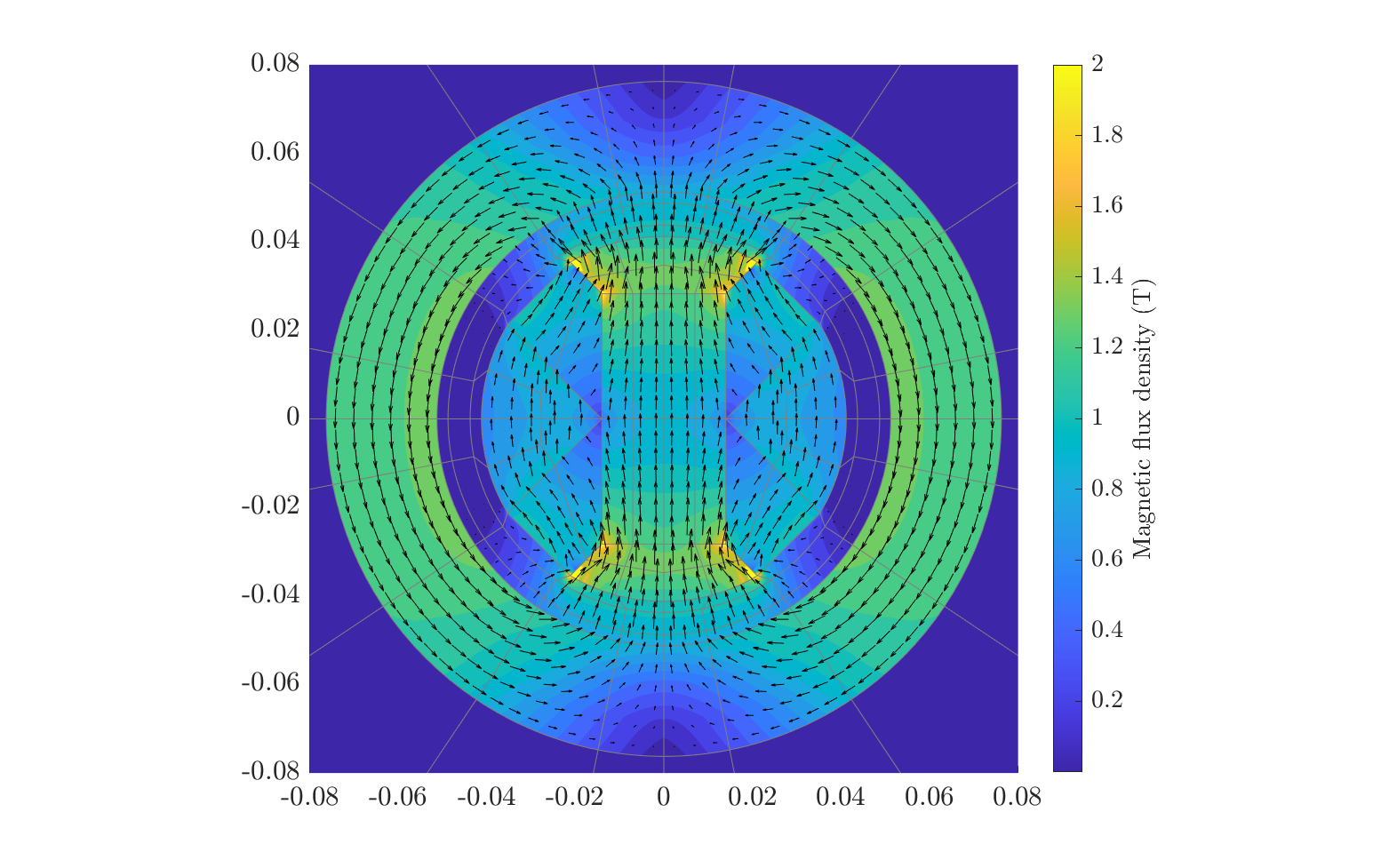}
      \caption{Magnetic flux density in Tesla.  For the final evaluation, the soft iron is implemented with a B-H-curve using fixpoint iteration.}
      \label{fig:optimization:resultsBField}
    \end{subfigure}
    \caption{Deduced geometry and analysis of the magnetic flux density.}
    \label{fig:optimization:results}
\end{figure*}

\Cref{fig:optimization:resultsGeometry} shows the simplified geometry created based on the topology optimization. The green patches indicate permanent magnets with the magnetization direction given by the arrows. The gray and blue patches indicate soft iron and air, respectively. The shapes of the permanent magnets are chosen in a way that allows an efficient production. In particular, the magnets can be easily cut from rectangular segments 
and magnetized uniformly, which is usually not the case for assemblies given in the literature \citep{BJORK2010437,MagnetocaloricEnergyConversion}.
The one remaining challenge with this design is to ensure the concentricity of the assembly in the absence of a continuous shaft.

\begin{figure}
    \centering
    \includesvg[width=\linewidth]{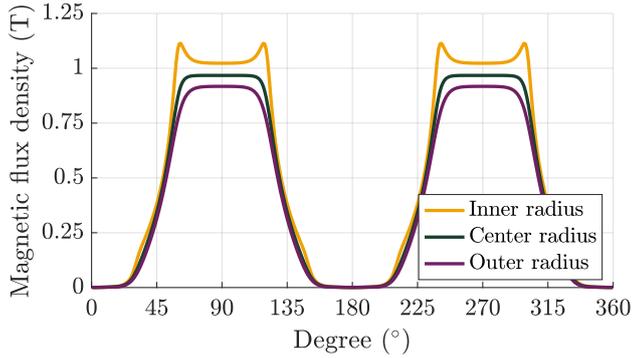}
\caption{Magnetic flux density in the empty air gap, evaluated at the lines indicated in \cref{fig:optimization:resultsGeometry}. }
    \label{fig:optimization:BGapTopOptGeometry}
\end{figure}

The resulting magnetic flux density is shown in \cref{fig:optimization:resultsBField}. With a share of  $A_\%=\SI{60}{\percent}$ of permanent magnets in the rotor, an objective value of $|\mathbf{B}|=\SI{0.76}{\tesla}$ is obtained, which is in very good agreement with the predicted value from the topology optimization.
\Cref{fig:optimization:BGapTopOptGeometry} shows the evaluation of the magnetic flux density in the air gap around the circles at $r = 42.5, 45$ and $\SI{47.5}{\mm}$ (see colored lines in \cref{fig:optimization:resultsGeometry}). At the center radius, the magnetic flux density is constant during magnetization with a value of $|\mathbf{B}|=\SI{0.96}{\tesla}$. In the radial direction, the magnetic flux density decreases as the flux is distributed across more space. At the inner radius, peaks in the magnetic flux density are observed. 
 
The magnetic flux density in the air gap is verified by nonlinear simulations carried out in COMSOL Multiphysics\textsuperscript{\textregistered}, which show maximum differences of less than \SI{0.1}{\percent}. Errors due to neglecting of saturation effects turn out to be smaller than \SI{0.5}{\percent}, which also justifies to carry out the optimizations with linear material properties.

Often, the coefficient 
\begin{equation}
    \Lambda _{\mathrm{cool}} =\left( \langle B_{\mathrm{high}}^{0.7} \rangle -\langle B_{\mathrm{low}}^{0.7} \rangle \right)\frac{V_{\mathrm{field}}}{V_{\mathrm{mag}}} P_{\mathrm{field}}
\end{equation}
is used as reference when comparing magnetocaloric cooling assemblies \citep{Bjrk2008}. It represents the ratio of magnetized field volume $V_\mathrm{field}$ to magnetizing volume $V_\mathrm{mag}$ multiplied by the average magnetic flux density in the high and low fields. The variable $P_\mathrm{field}$ describes the fraction of a cooling cycle where MCM is placed in the high field. The magnetic flux density exponent of 0.7 comes from the observation that the temperature change of Gadolinium at the Curie temperature is approximately proportional to $B^{0.7}$. 

In our case the the magnetic field densities are
\begin{equation}
\langle B_\mathrm{high}^{0.7} \rangle = \SI{0.97}{\tesla}
\quad
\text{and}
\quad
\langle B_\mathrm{low}^{0.7} \rangle = \SI{0.03}{\tesla}.
\end{equation}
In the presented case of 2D simulations, the volume ratio equals the ratio of the 2D areas. 
With a magnet area of \SI{31.9}{\centi\metre\squared}, an air gap area of \SI{29.8}{\centi\metre\squared} of which one third is magnetized,
$P_\mathrm{field} = 0.75$, it follows that $\Lambda_\mathrm{cool} = 0.219$ and $\Lambda_\mathrm{cool}/P_\mathrm{field} = 0.29$. These values are even higher compared to the ones provided in the literature \citep{BJORK2010437}. The highest value of $\Lambda_\mathrm{cool} = 0.21$ was reported from \cite{Okamura}. 
Here, these values are given only for comparison as $\Lambda_\mathrm{cool}$ has some drawbacks, e.g., it does not consider MCM cost, cooling power or the device efficiency \citep{benke2021}.
An even higher $\Lambda_\mathrm{cool}$ could be achieved here by increasing the air gap.
It has therefore been proposed to include $B^{4/3}$ in such merit figures~\citep{benke2020}.

\subsection{Shape Optimization}
Measurements on a magnetocaloric device indicate that peaks in the magnetic flux density during the magnetization reduce the performance of the MCM \citep{PerformanceMeasurementsOnALargeScale}. A lower but more homogeneous magnetic field performed significantly better than a higher,  inhomogeneous one. In \citep{PerformanceMeasurementsOnALargeScale}, the authors stated that these peaks create ``parasitic thermodynamic cycles'' and ``represent additional irreversibilities''. Also in \cref{fig:optimization:BGapTopOptGeometry}, we observe peaks during magnetization at the inner evaluation radius. To address this issue at an early design stage, shape optimization is applied to the soft iron parts of the geometry in order to remove the peaks in the magnetic flux density.

\subsubsection{Methodology}
\label{sec:ShapeOptimization:Methodology}

So far, the influence of the MCM has widely been neglected, as the standard material for scientific prototypes, Gadolinium, reaches saturation very quickly and behaves similarly to air regarding the resulting magnetic field. As this device is aimed for commercialization, Gadolinium is not an option due to its high cost. 
Current MCM candidates for commercialization are alloys based on Lanthanum, Iron and Silicium (La-Fe-Si), which have a higher permeability. Thus, neglecting the MCM's influence is not an option.

While some measurements of magnetization curves exist \citep{Gozdur2020}, the behavior of the MCM depends on many factors, such as the ambient and Curie temperatures. We therefore use a worst case estimate for the permeability to investigate this problem.
Assuming a conservative relative permeability of $\mu_\mathrm{r, MCM}=\SI{500}{}$ and a volume fraction of \SI{60}{\percent} for the MCM with water as second phase with $\mu_\mathrm{r, water}=1$, the upper limit of the permeability is given by the Voigt bound \eqref{eq:modeling:multiscale:voigt2} as $\mu_{\mathrm{Voigt}} = 300$, which corresponds to parallel plates. For packed spheres, $\mu_\mathrm{r}$ is consequently smaller, so we use an isotropic $\mu_\mathrm{r, AMR} = 200$ for the AMRs to simulate the edge cases of the macroscopic behavior.

The aim is again to maximize $\langle B_{\mathrm{high}} \rangle  - \langle B_{\mathrm{low}} \rangle$, i.e., the difference in magnetic flux density between the magnetization and the demagnetization phase. 
The problem is defined in a similar way to the topology optimization as
    \begin{align}
    \max_{\mathbf{P}} \quad & f_{\opt} = \langle B_{\mathrm{high}} \rangle  - \langle B_{\mathrm{low}} \rangle \\ 
    \text{s.t.} \quad & \norm{\mathbf{P}_\mathrm{i}}   \leq  R_{\mathrm{max}}\\
      & \mathbf{\| B}_\mathrm{eval}(\beta_{i}) \| \leq \mathbf{\| B}_\mathrm{eval}(\beta_{i-1}) \|     \\
      \text{and} \quad &  \text{magnetostatics~\eqref{eq:NumericalMethods:Mortaring:MatrixSystem} for all $\beta_i$},
    \end{align}
where this time the design variables are the control points $\mathbf{P}_\mathrm{i}$ of the geometry parametrization described in \cref{sec:geometry}.  Instead of evaluating multiple points in the air gap, one point in the MCM is evaluated for different rotation angles $\beta \in \{0,1,\dots 90\}$.\footnote{The optimization has also been performed with an empty air gap, but since the results are very similar, only the results for a filled air gap are presented.} The values $\langle B_{\mathrm{high}} \rangle$ and $\langle B_{\mathrm{low}} \rangle $ are then obtained by averaging the magnetic flux density for the rotation angles during magnetization and demagnetization, respectively.

Again, constraints are needed. First, the control points are restricted by the radius $R_\mathrm{max}$ to ensure feasible designs. 
Second, the magnetic flux density at each previous evaluation angle $\beta_{i-1}$ must be larger than the magnetic flux density at the current evaluation angle $\beta_{i}$. This is sufficient, because only one (de)magnetization process is simulated due to the symmetries. Consequently, the magnetic flux density must be monotonically decreasing and become more homogeneous, as peaks are infeasible.

\begin{figure}
    \centering
    \includegraphics[trim= 210 125 184 100, clip,
    width=\linewidth]{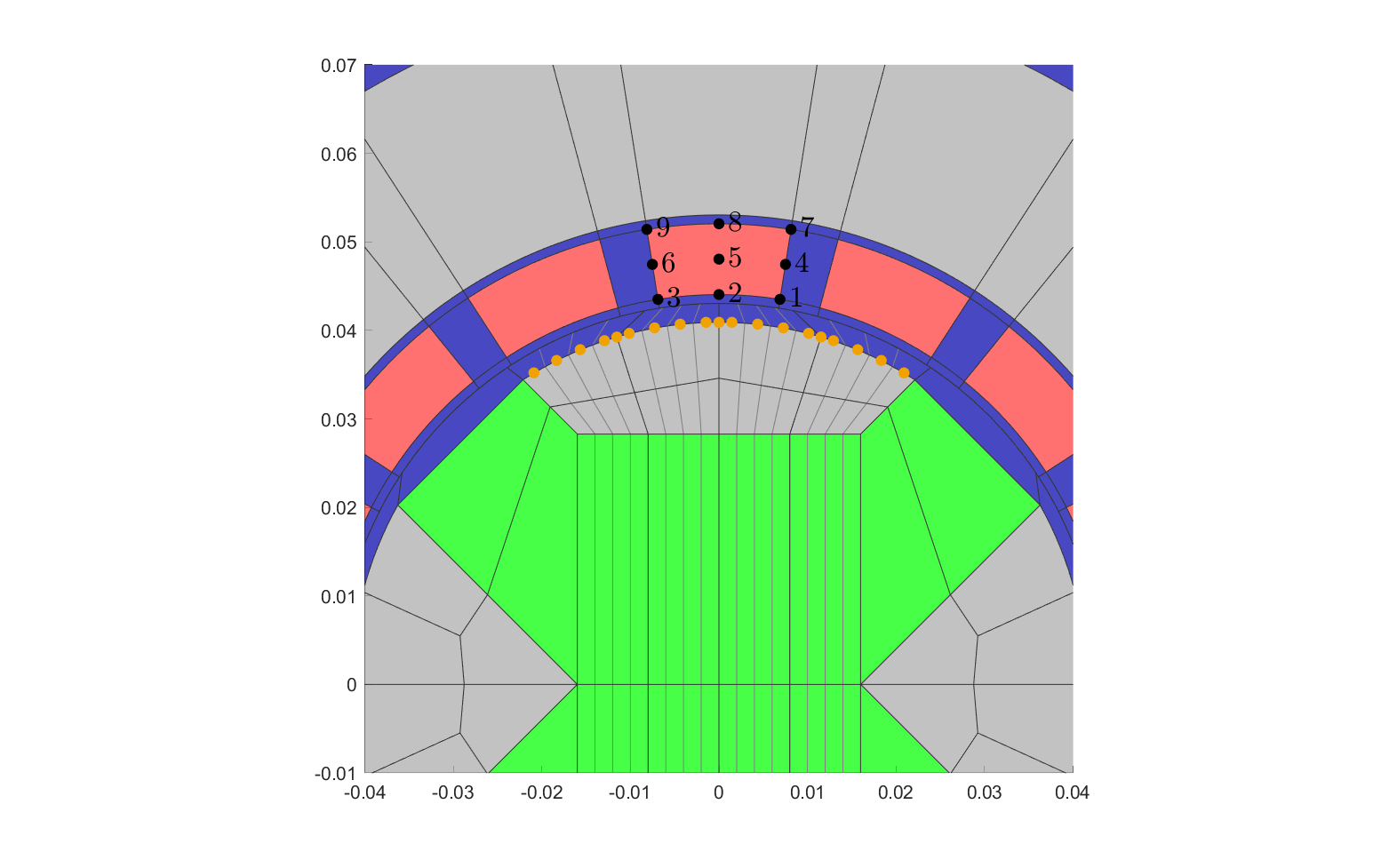}
    \caption{Initial geometry for the shape optimization. 
    The orange points denote the control points, which serve as design variables. 
    }
    \label{fig:Optimization:ShapeOpti:Initial}
\end{figure}

The initial geometry is visualized in \cref{fig:Optimization:ShapeOpti:Initial}.  The control points, which define the soft iron surface and act as design variables for the shape optimization, are shown in orange. Soft iron and magnets are again shown as gray and green patches, where the magnetization direction is the same as in \cref{fig:optimization:resultsGeometry}. Blue patches indicate air or any support material behaving like air.
The magnetic flux density is evaluated at the black points in the AMR, which are the red patches. Point 2 is used to evaluate the magnetic flux density for the rotation angles from $\beta = \SI{0}{\degree}$ to \SI{90}{\degree} in \SI{1}{\degree} steps.  This simulates the process of the top center AMR from being fully magnetized to demagnetized.

The optimization is again realized with the MATLAB\textsuperscript{\textregistered} solver \textit{fmincon} using analytically derived gradients for the optimization function and constraints. For the shape optimization, the  \textit{SQP} method is used, as it converges faster for a low number of design variables and is better suited for infeasible points, since the initial condition violates the second constraint due to the peak. 

\subsubsection{Results}

\begin{figure}
    \centering
    \includesvg[width=\linewidth]{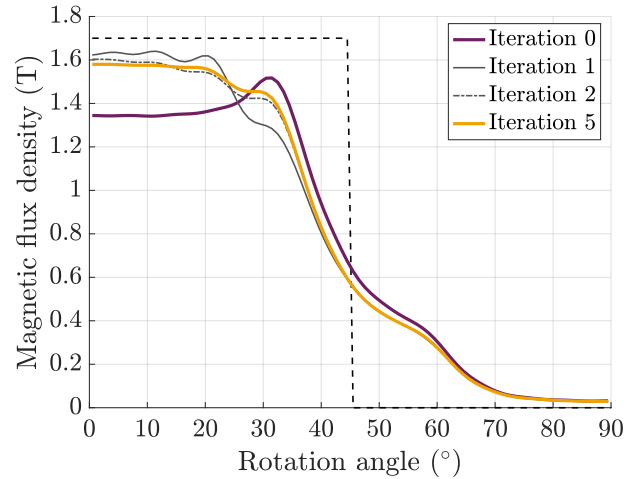}
    \caption{Magnetic flux density in the AMR over the course of the optimization. Starting with a peak in the initial design, the magnetic flux density is gradually homogenized during the optimization process, resulting in a higher and more uniform magnetic flux density. }
    \label{fig:Optimization:ShapeBValues}
\end{figure}

The magnetic flux density during the optimization process at the evaluation point is shown in \cref{fig:Optimization:ShapeBValues}.  In the initial design, the magnetic flux is concentrated in the AMR part which has already entered the magnetization zone. As the whole AMR enters the magnetization zone, the magnetic flux density is distributed across the entire AMR volume. This effect results in a falling magnetic flux density in the AMR, causing the peaks during magnetization. 
For the optimized shape, the distance of the AMR to the soft iron changes during operation. The smaller distance at the center increases the magnetic flux density there and shifts the high concentration of the magnetic flux density to the center. The peaks during magnetization are therefore removed, leading to a more homogeneous magnetic field in the AMR. Note that the magnetic flux density in the AMR is simulated, resulting in a higher magnetic flux density in \cref{fig:Optimization:ShapeBValues} compared to \cref{fig:optimization:BGapTopOptGeometry} due to the higher $\mu_\mathrm{r}$.

\begin{figure*}
    \centering
    \begin{subfigure}[t]{.32\linewidth}
      \centering
        \includegraphics[trim= 250 225 224 50, clip, width=\linewidth]{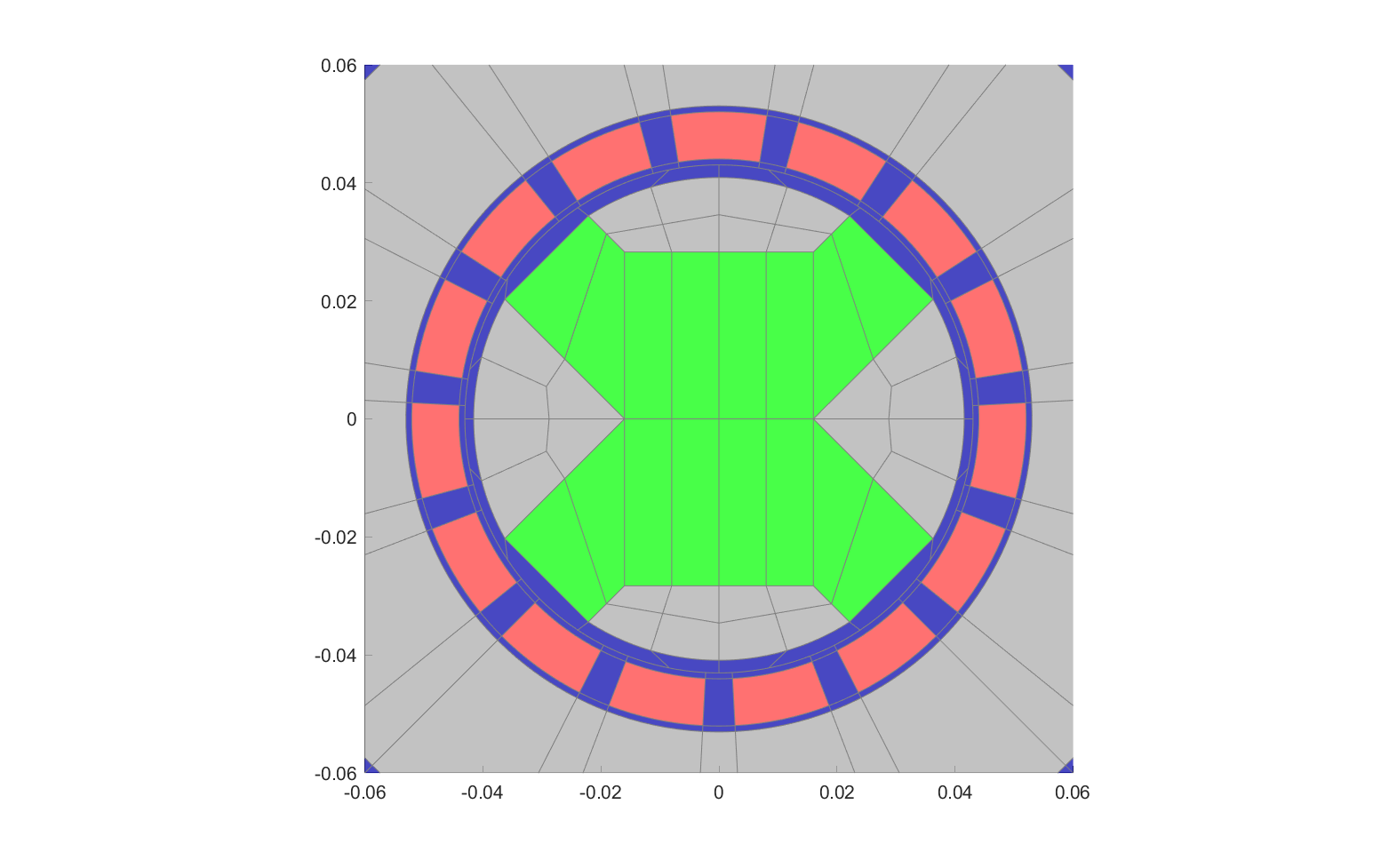}
        \caption{Iteration 0.}
      \label{fig:optimization:ShapeOpt:Results1}
    \end{subfigure}
    \begin{subfigure}[t]{.32\linewidth}
      \centering
        \includegraphics[trim= 250 225 224 50, clip, width=\linewidth]{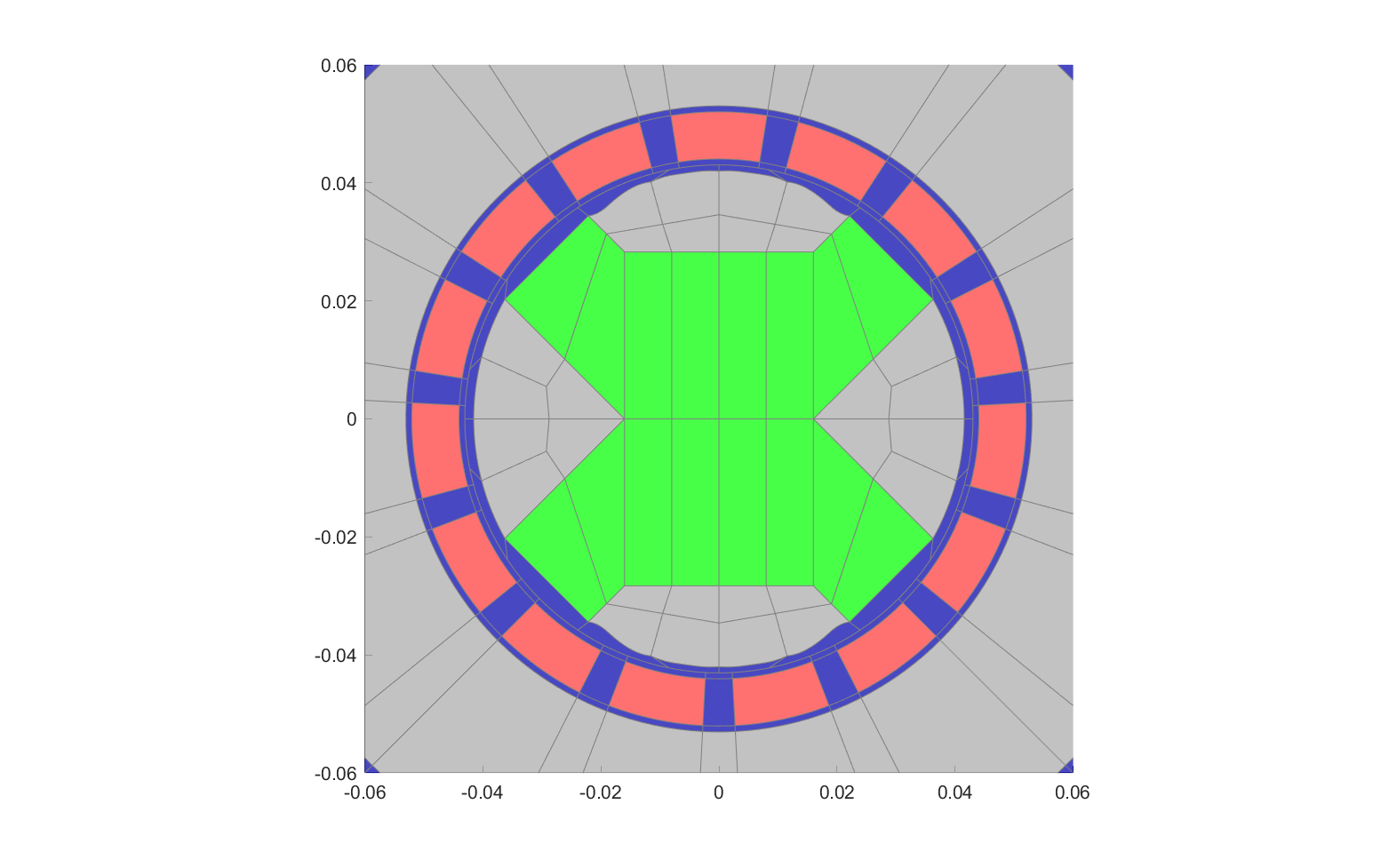}
      \caption{Iteration 1.}
      \label{fig:optimization:ShapeOpt:Results2}
    \end{subfigure}
    \begin{subfigure}[t]{.32\linewidth}
      \centering
        \includegraphics[trim= 250 225 224 50, clip, width=\linewidth]{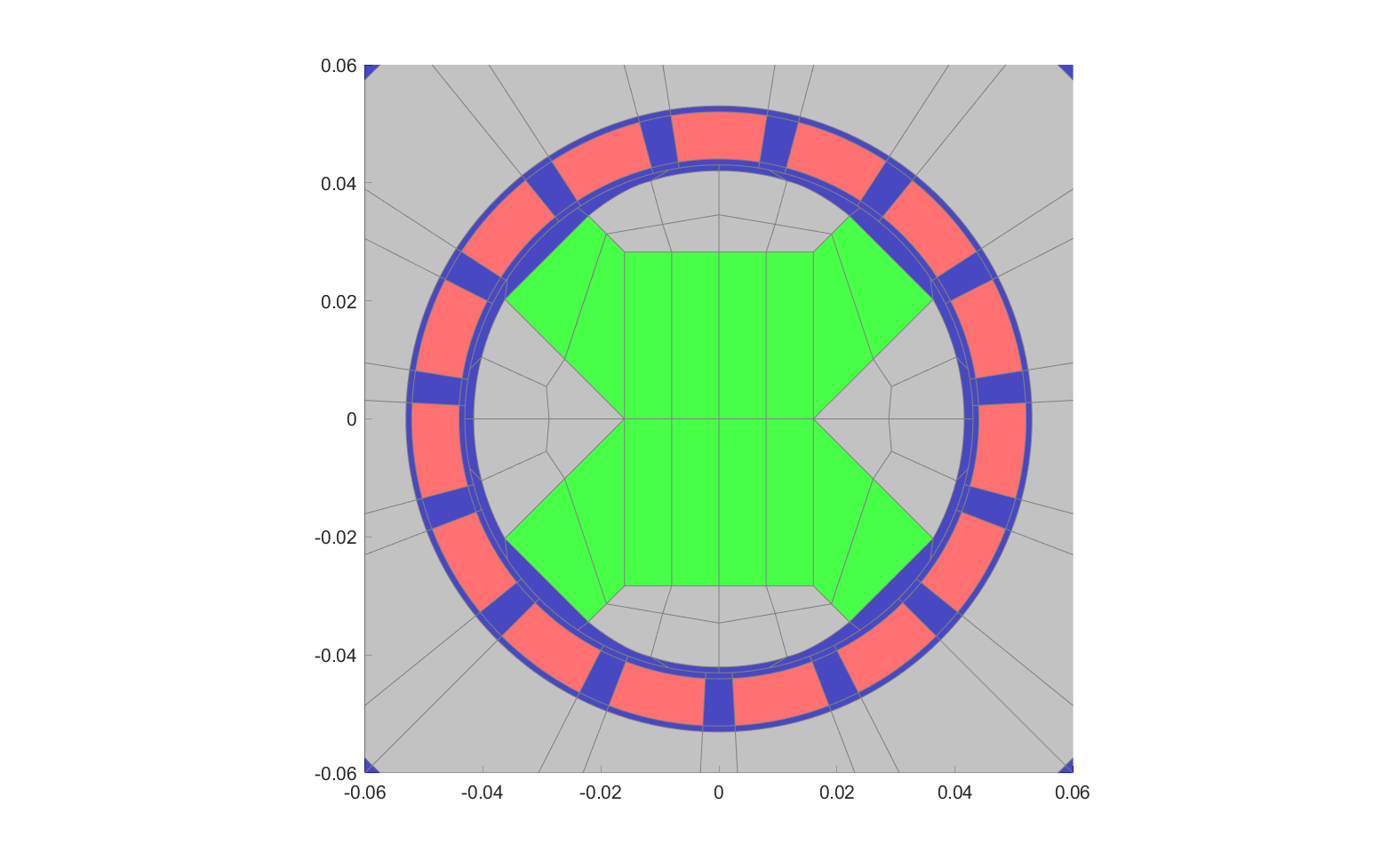}
      \caption{Iteration 5.}
      \label{fig:optimization:ShapeOpt:Results3}
    \end{subfigure}
    
    \caption{Shape optimization process for a filled air gap with simplified AMR properties. A circular arc with center at the coordinate origin is used as initial condition.
    In the first iteration, the solver increases the overall radius, except for the peak areas, where the AMR enters (or leaves) magnetization. There, the radius is decreased to remove the peaks in the magnetic flux density. After several iterations, the optimizer converges to a design, where the distance to the AMR is gradually decreased as the system rotates.}
    \label{fig:optimization:ShapeOpt:results}
\end{figure*}

The iteration process for the shape optimization is shown in \cref{fig:optimization:ShapeOpt:results}. Starting from a circular arc with its center at the coordinate origin as initial condition (\cref{fig:optimization:ShapeOpt:Results1}), the radius of the soft iron is increased in the center to yield a higher objective function and decreased at the left and right sides to avoid the peaks in the magnetic flux density (\cref{fig:optimization:ShapeOpt:Results2}). The solver quickly converges to a design where the distance of the AMR to the soft iron is slowly reduced (\cref{fig:optimization:ShapeOpt:Results3}), yielding a smooth transition of the AMR into the magnetization zone when the system rotates. Increasing the distance to the AMR at the magnet edges has also been shown to increase the magnetic flux density in other assemblies \citep{cryst9020076}.

For a practical design, it is possible to approximate the spline geometry of the rotor well by a circular arc where the circle center is not the origin, but shifted in the $y$-direction. This may simplify the manufacturing process.

\section{Final design studies}
\label{chap:FinalDesignStudies}

\subsection{Torque study}

For an efficient refrigeration system, it is essential to minimize  occurring torques during operation to reduce the power consumption of the motor \citep{Bouchekara2011}. Also, the rotation becomes more stable and valve mistimings of the hydraulic system are less likely to occur.

As the magnetic flux is concentrated in the AMRs due to the increased permeability, the magnetic flux is not evenly distributed in the air gap, resulting in reluctance forces and torque. This depends mainly on the geometric configuration of the rotor and stator, i.e., the setup of the permanent magnets and the AMRs \citep{LOZANO2016187}. 
This is particularly important for devices aimed to be operated with La-Fe-Si, as the higher permeability leads to torques several times higher compared to Gadolinium. We only examine the static holding torque (cogging torque) and neglect effects like eddy currents or friction. 

\begin{figure}
    \centering
   \includesvg[width=\linewidth]{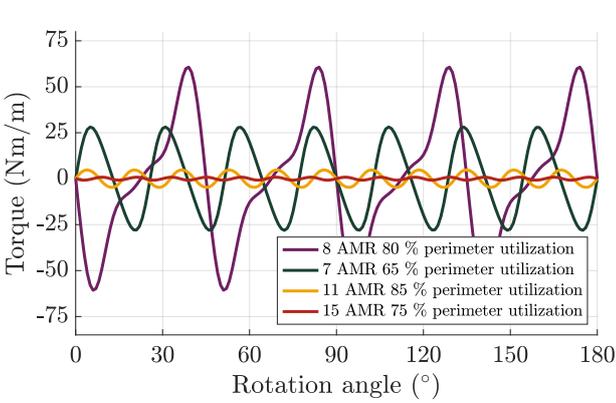}
    \caption{Exemplary torque curves for different stator setups, varying the number of AMRs and the share of MCM along the perimeter in the air gap. }
    \label{fig:NumericalResults:Torques}
\end{figure}

The rotor geometry has already been determined with the topology optimization, so a parametrized stator geometry is analyzed with the initial rotor design from \cref{fig:Optimization:ShapeOpti:Initial}. 
There are two parameters that have a major influence on the cogging torque, that is, the number of AMRs $\NAMR$ and the share of AMRs along the perimeter $\rAMR$ (called perimeter utilization in the following). The AMRs are designed to be evenly distributed along the perimeter.   \Cref{fig:NumericalResults:Torques} shows exemplary torque curves for half a revolution. The values are calculated per unit length, so the length of the device needs to be multiplied for the real values. For an even number of AMRs, the number of maxima and minima corresponds to the number of AMRs. For an odd number of AMRs, there are twice as many maxima and minima, since the magnetic field at the poles is not symmetric. In this case, opposing sides do not necessarily contribute forces in opposite directions, resulting in considerably smaller torques. 
During the rotation, when one AMR is moving out of the magnetization zone, it creates a torque counteracting the rotation direction. In a scenario with an even $\NAMR$, due to symmetry, the same is true for the opposite side of the rotor such that the torque is summed up. However, if $\NAMR$ is odd, the AMR on the opposite side instead enters the other magnetization zone. This results in a torque in the different direction and the total torque is reduced. The decreasing torque values for an increasing $\NAMR$ are explained by their smaller areas which reduces the magnetic flux per AMR and therefore their associated reluctance torque.

\begin{figure}
    \centering
    \includegraphics[width=\linewidth]{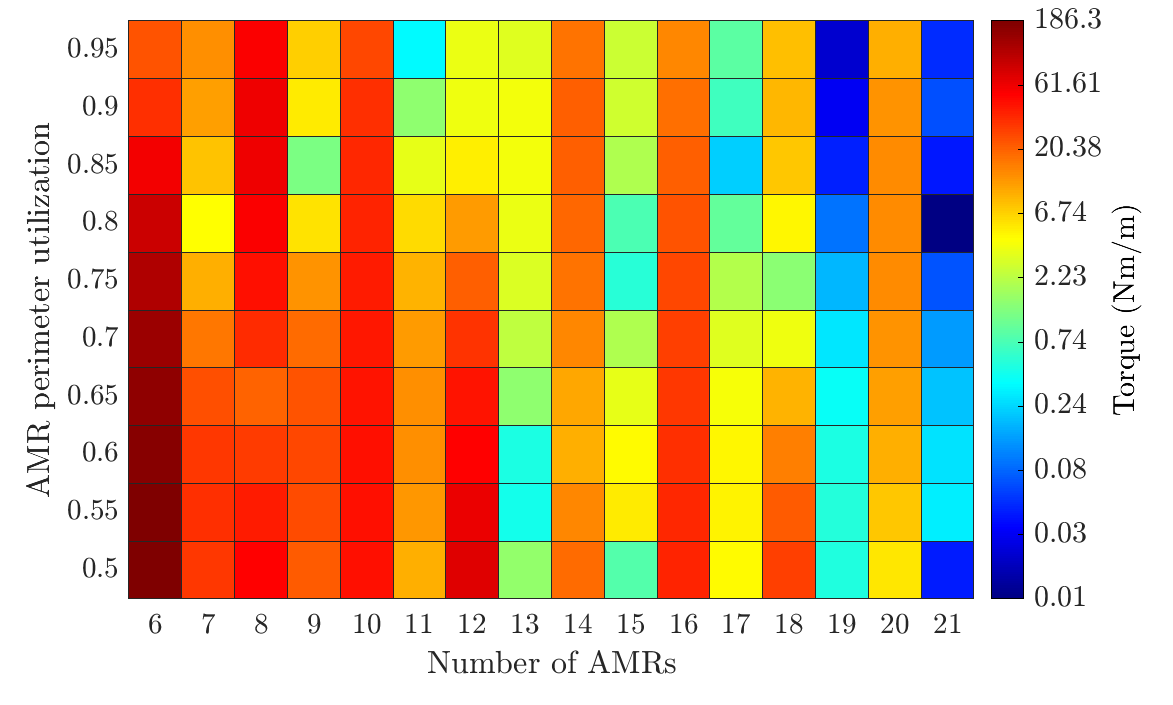}
    \caption{Holding torque per system length, depending on the number of AMRs and the perimeter utilization. Higher and odd numbers of AMRs correspond to lower torques.}
    \label{fig:NumericalResults:TorquesHeatmap}
\end{figure}

Performing the torque simulations for a variety of parameter combinations yields the heat map presented in \cref{fig:NumericalResults:TorquesHeatmap}.
The following observations are made:
\begin{itemize}
    \item The most influential parameter is the number of AMRs $\NAMR$. Setups with an odd $\NAMR$ have significantly lower torques than setups with an even $\NAMR$. A higher $\NAMR$ also leads to lower torques.
    \item For a fixed $\NAMR$, there are again torque minima depending on $\rAMR$. They are explained by geometric setups, where one AMR enters the magnetization zone exactly when another one leaves.
\end{itemize}
Thus, in terms of torque, it is advisable to have an odd, high number of AMRs. Unless the MCM mass is not already determined by the specifications of the device, $\rAMR$ can be further adjusted to reduce the torque. 
This allows for a minimization of the cogging torque, which is more problematic with, e.g., Halbach cylinder designs \citep{Arnold2014}.

\subsection{Magnetic field evaluations}

In the following, we present three different examples to highlight positive and negative design aspects. Note that these scenarios pose edge cases where the MCM  permeability is overestimated to demonstrate its influence and to deduce general design rules.
\begin{figure*}[t]
    \centering
    \begin{subfigure}[t]{.32\linewidth}
      \centering
        \includegraphics[trim= 220 102 252 88, clip, width=\linewidth]{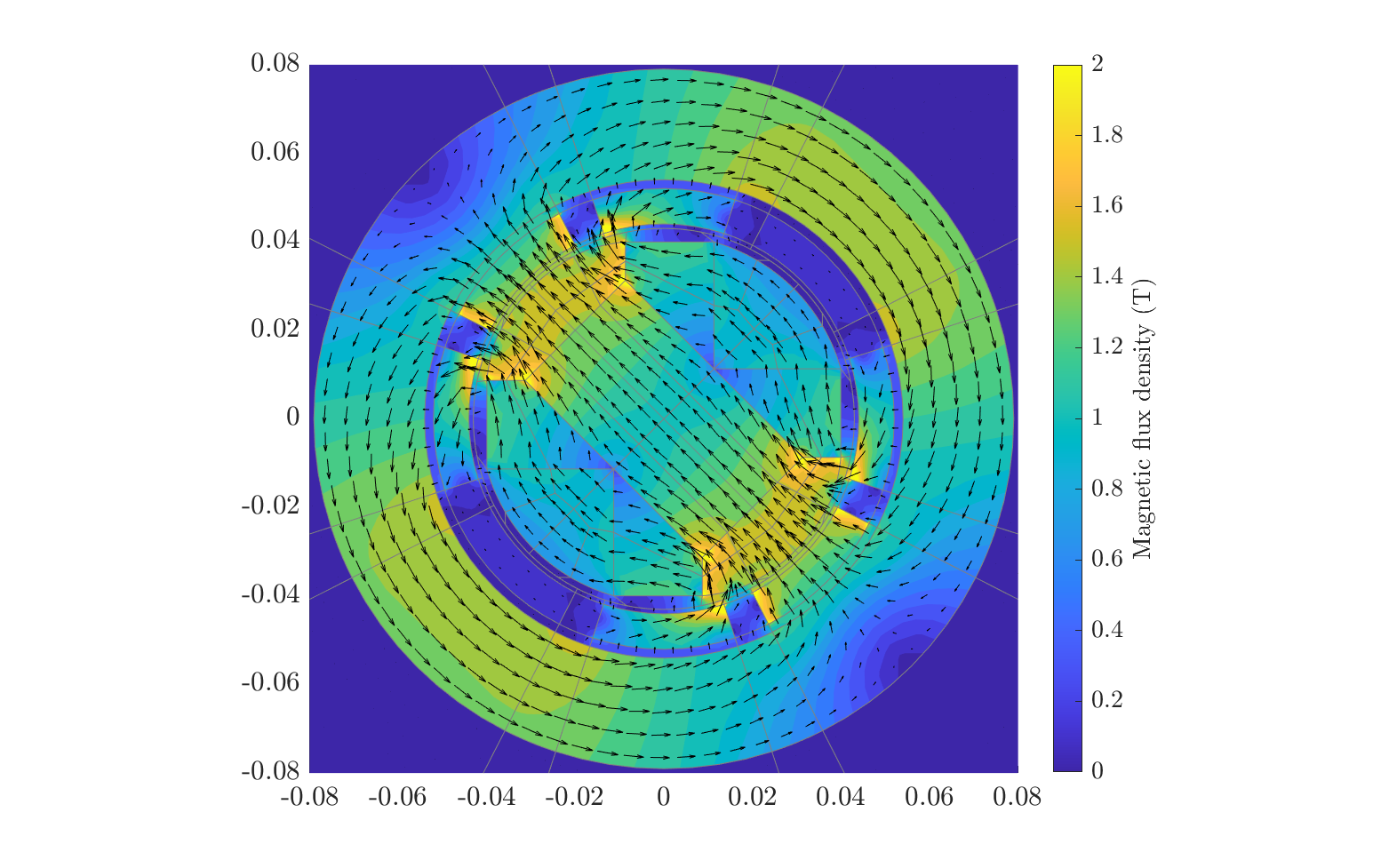}
        \caption{8 isotropic AMRs with $\rAMR = \SI{80}{\percent}$. }
      \label{fig:NumericalResults:8AMR}
    \end{subfigure}
    \begin{subfigure}[t]{.32\linewidth}
      \centering
        \includegraphics[trim= 220 102 252 88, clip, width=\linewidth]{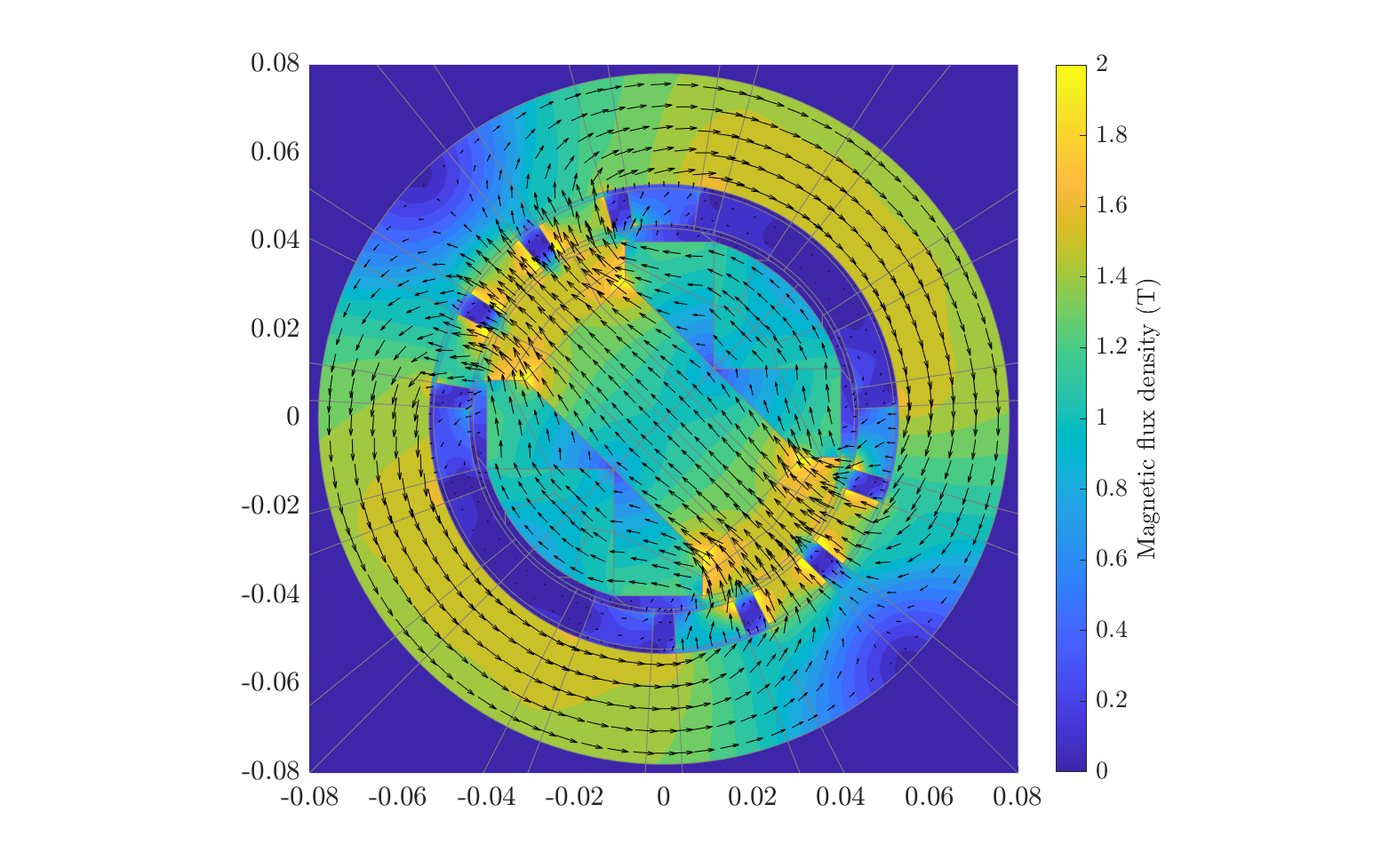}
      \caption{15 isotropic AMRs with $\rAMR = \SI{75}{\percent}$.}
      \label{fig:NumericalResults:15AMR}
    \end{subfigure}
    \begin{subfigure}[t]{.32\linewidth}
      \centering
        \includegraphics[trim= 220 102 252 88, clip, width=\linewidth]{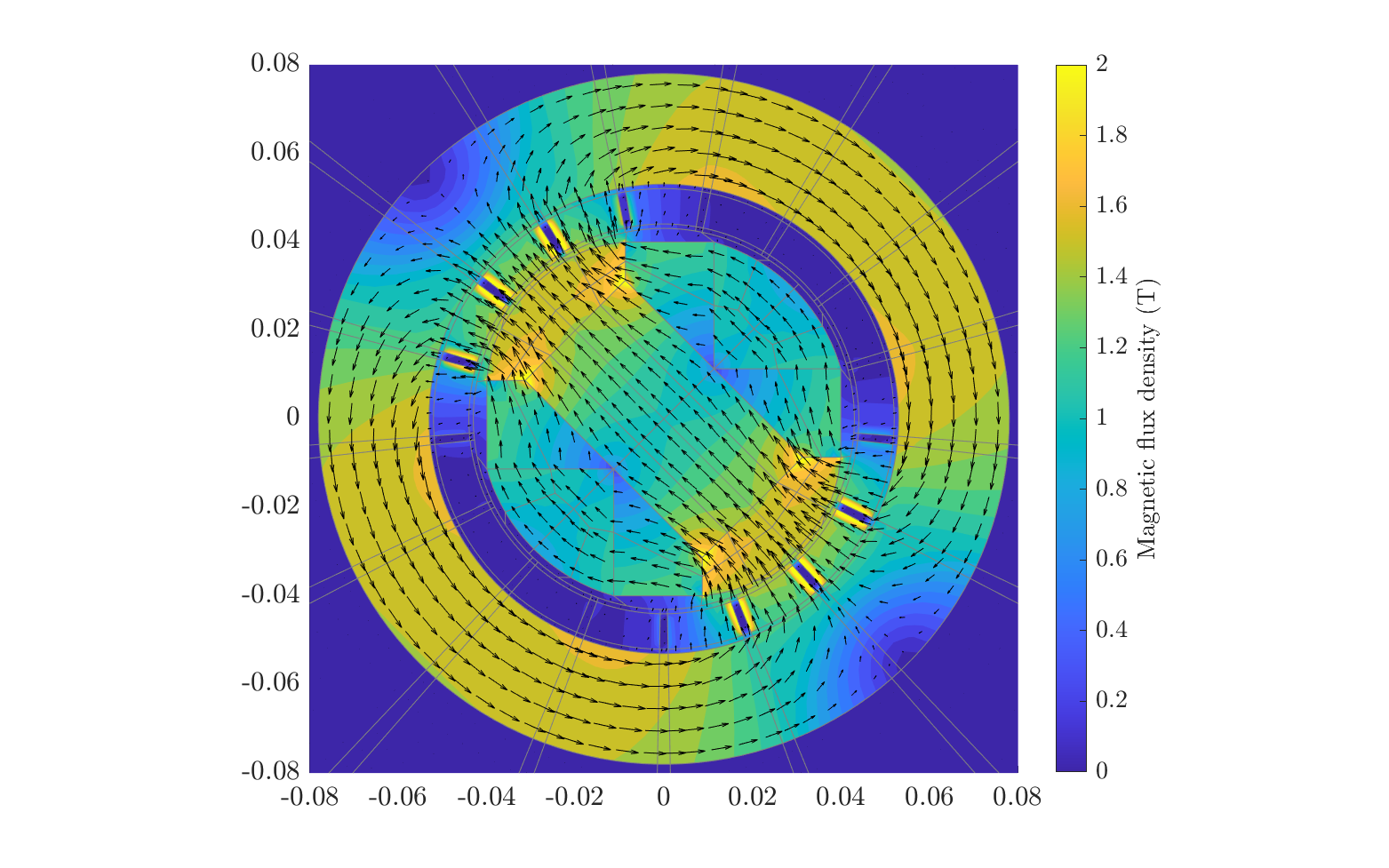}
      \caption{17 anisotropic AMRs with $\rAMR = \SI{90}{\percent}$.}
      \label{fig:NumericalResults:17AMR}
    \end{subfigure}
    
    \caption{Magnetic flux density for three exemplary AMR designs. (a) An unfavorable design where the AMRs are too long, resulting in their magnetization, despite being only partially in the magnetization zone. (b) A favorable design with a higher $\NAMR$ and lower $\rAMR$ resulting in more distinct (de)magnetization zones. (c) A prospective scenario with both a high $\NAMR$ and $\rAMR$. The magnetic flux is radially enforced due to the anisotropy of the AMRs, which allows distinct (de)magnetization zones and a high usage of MCM in the air gap.
    }
    \label{fig:NumericalResults:BFields}
\end{figure*}

\begin{figure*}
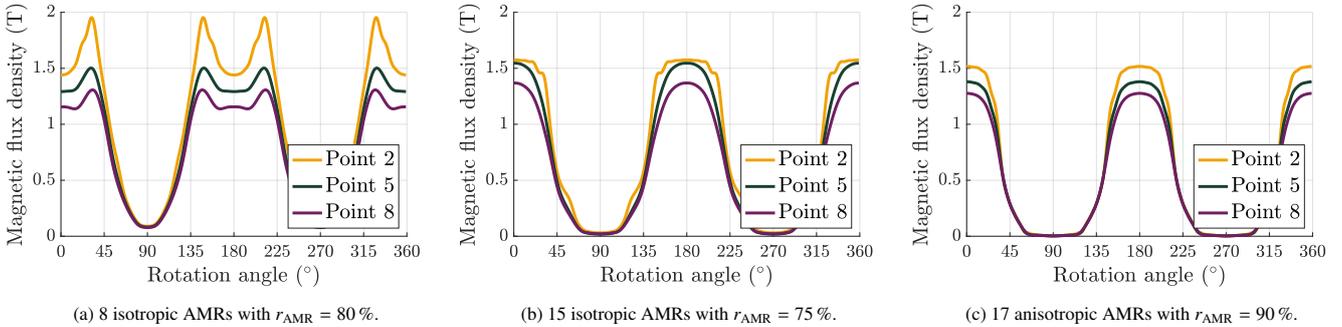

    \centering
    \begin{subfigure}[t]{.32\linewidth}
      \centering
        \includesvg[width=\linewidth]{8AMRPoints258}
        \caption{8 isotropic AMRs with $\rAMR = \SI{80}{\percent}$.}
      \label{fig:NumericalResults:8AMRFlux}
    \end{subfigure}
    \begin{subfigure}[t]{.32\linewidth}
      \centering
        \includesvg[width=\linewidth]{15AMRPoints258}
      \caption{15 isotropic AMRs with $\rAMR = \SI{75}{\percent}$.}
      \label{fig:NumericalResults:15AMRFlux}
    \end{subfigure}
    \begin{subfigure}[t]{.32\linewidth}
      \centering
        \includesvg[width=\linewidth]{17AMRPoints258}
      \caption{17 anisotropic AMRs with $\rAMR = \SI{90}{\percent}$.}
      \label{fig:NumericalResults:17AMRFlux}
    \end{subfigure}
    
    \caption{Magnetic flux density in the AMRs evaluated at the points given in \cref{fig:Optimization:ShapeOpti:Initial}. (a) An unfavorable design with high peaks in the magnetic flux density during the magnetization and a non-zero magnetic flux density during demagnetization. (b) A favorable design with a very homogeneous magnetic flux density during magnetization and a zero-field during demagnetization. (c) The highest MCM usage in the air gap and the most favorable magnetic flux density, which is reached with anisotropic AMRs.}
    \label{fig:NumericalResults:Flux}
\end{figure*}

The first example features a design with $\NAMR=8$ AMRs and $\rAMR = \SI{80}{\percent}$, where the AMRs are assumed to be isotropic with a constant relative permeability of $\mu_\mathrm{r, AMR}=200$, as explained in \cref{sec:ShapeOptimization:Methodology}. A lower number of AMRs is sensible from a systemic point of view in order to simplify the hydraulic system. However, looking at  the magnetic flux density in \cref{fig:NumericalResults:BFields}, we see that this design is not preferable from a magnetic field perspective for multiple reasons. First, the AMRs are too wide. This results in magnetization despite only being partially in the magnetization zone. Second, the high $\rAMR$ leads to a parasitic backflow of the magnetic flux through the AMRs. This results in a very narrow demagnetization zone, where the magnetic flux density is not zero.
Third, the magnetic flux is concentrated when the AMR enters the magnetization zone, which leads to high peaks. For wide AMRs, this cannot be resolved with shape optimization, which is why the optimized rotor for $\NAMR=15$ AMRs is used for the evaluations. Last, the torque is very high due to the symmetric setup, as seen in \cref{fig:NumericalResults:Torques,fig:NumericalResults:TorquesHeatmap}.

A better approach is to use more AMRs and a lower perimeter utilization, as performed for the second design shown in \cref{fig:NumericalResults:15AMR} for $\NAMR=15$ and $\rAMR=\SI{75}{\percent}$. This ensures that the magnetic flux density in the AMRs is zero during demagnetization and there are no peaks during magnetization, as seen in  \cref{fig:NumericalResults:15AMRFlux}. As a result, long (de)magnetization phases and a homogeneous flux density are achieved.
In general, it is advisable to choose the width of the AMRs as well as the thickness of the air gap smaller than the width of the diagonal magnets to ensure the magnetic flux flow through the stator iron. Also note that the magnetic flux density in point 2 is exactly the final curve from \cref{fig:Optimization:ShapeBValues}, since it was optimized for this case. The torque is also very small, as seen again in \cref{fig:NumericalResults:Torques,fig:NumericalResults:TorquesHeatmap}. Another advantage of having more and smaller AMRs is that the shape can be simplified to a rectangle, which might be convenient for AMRs with ordered microstructures.

This case of anisotropic material behavior for the AMRs is investigated in the third example with $\NAMR=17$ and $\rAMR=\SI{90}{\percent}$. Assuming radially oriented phases, the radial and azimuthal permeability is calculated with (\ref{eq:modeling:multiscale:voigt2}) and (\ref{eq:modeling:multiscale:reuss2}) to $\mu_{\mathrm{r, AMR, rad}}=300$ and $\mu_{\mathrm{r, AMR, azi}} = 2.5$ under the previous assumptions of $\mu_\mathrm{r, MCM} = \SI{500}{}$ and a volume fraction of \SI{60}{\percent} for the MCM material. Generating this anisotropic behavior in reality could be realized with radially aligned plates or pins \citep{PerformanceAssessmentOfDifferentPorousMatrix}.
The permeability tensor for the AMRs is then given by (\ref{eq:modeling:MuTensor}). 
The resulting magnetic field is shown in \cref{fig:NumericalResults:17AMR}. Despite having a high $\rAMR$, there is no magnetic flux during the demagnetization, as azimuthal flux is impeded by the microstructure.
Clever exploitation of the MCM anisotropy leads to the fastest transition from the (de)magnetization phases and a homogeneous magnetic flux density, seen in \cref{fig:NumericalResults:17AMRFlux}, regardless of whether the rotor is shape-optimized.

On the other hand, this analysis implies to avoid microstructures that obstruct the magnetic flux in radial direction. Measurements of parallel plated Gadolinium AMRs, where the cooling performance is found to be up to \SI{40}{\percent} smaller for plates perpendicularly aligned to the magnetic field, support this theory \citep{TUSEK201357}. This effect is assumed to be even stronger for materials with higher permeability such as La-Fe-Si, because the magnetic field is concentrated even more in phases with higher permeability. Using microstructures correctly therefore not only positively influences the magnetic field, but also concentrates the magnetic flux in the MCM phases, which enhances the temperature generation in the AMRs.

\subsection{3D simulations}

\begin{figure}[t]
    \centering
    \includegraphics[width=\linewidth]{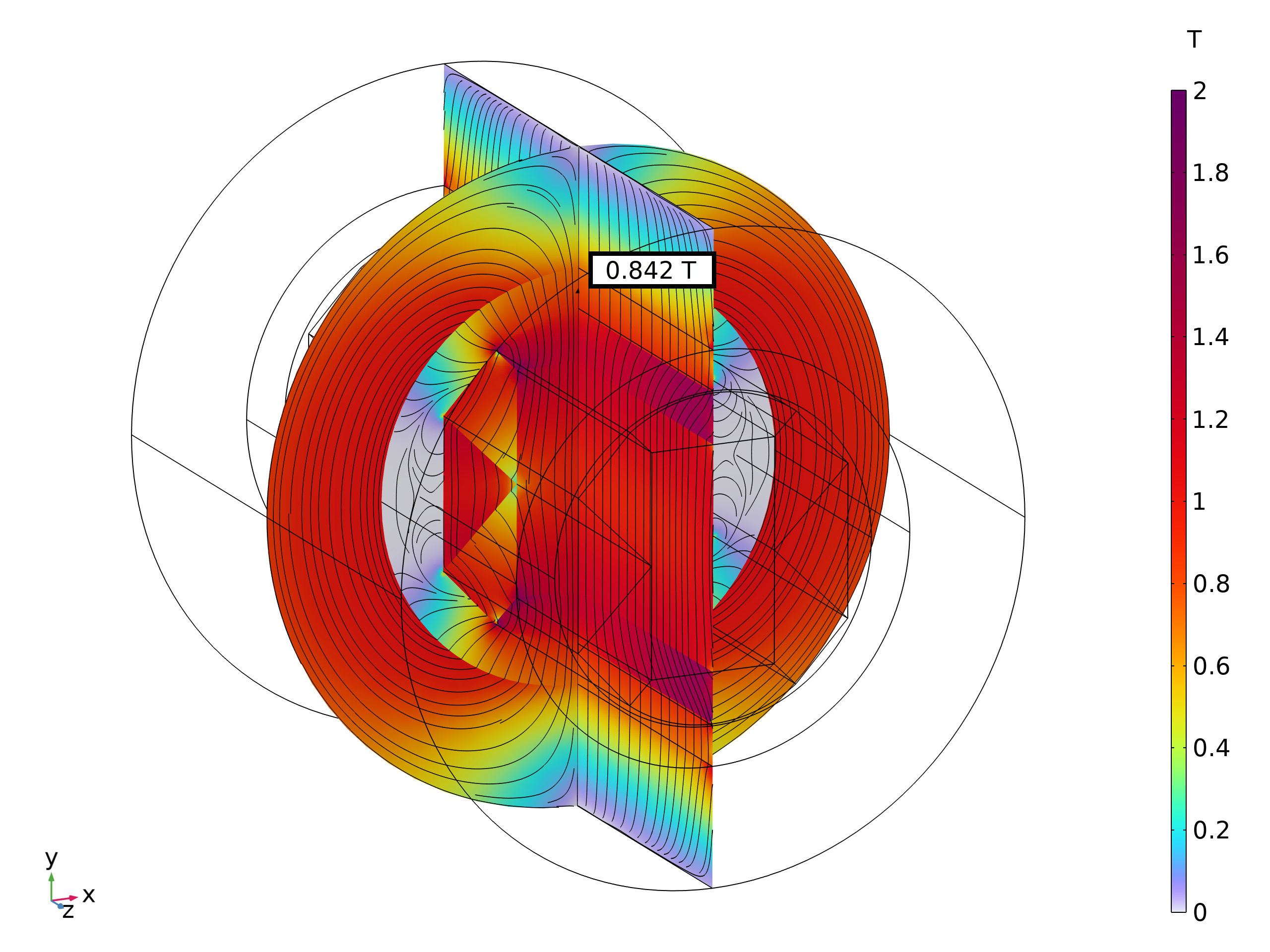}
    \caption{Magnetic flux density for the 3D case from \cref{fig:optimization:ShapeOpt:Results3} with an empty air gap, simulated with COMSOL Multiphysics\textsuperscript{\textregistered}.}
    \label{fig:Evaluation:3Dsimulation}
\end{figure}

It should be noted that 2D simulations usually overestimate the magnitude of the magnetic flux density as 3D fringe effects are neglected. 
A 3D simulation for the optimized machine from \cref{fig:optimization:ShapeOpt:Results3} with an empty air gap has been carried out in COMSOL Multiphysics\textsuperscript{\textregistered} to obtain a more realistic result and assess the differences caused by these fringe effects. Given an exemplary length of \SI{15}{\cm}, the mean magnetic flux density during magnetization is found to decrease from \SI{0.92}{\tesla} to an average of \SI{0.83}{\tesla} in the high field, which is about \SI{10}{\percent} smaller compared to the 2D simulation. The results from the 3D simulation are found in \cref{fig:Evaluation:3Dsimulation}.

Another simplification commonly made for the 2D-case is the neglecting of eddy currents as sufficiently small lamination is assumed. Here, a lamination thickness of \SI{0.5}{\mm}, as used in electric motors, is recommended to reduce the heat generation in the stator due to eddy currents. Although the frequency of the magnetic field at around \SIrange{5}{10}{\Hz} is much lower than that of electric machines, the heat generation in the stator is more harmful. This is because it enters the cooling circuit as a parasitic heat input, making a small lamination thickness necessary.

\section{Conclusion}
\label{chap:Conclusion}

This paper introduces advanced simulation and optimization methods to find the most suitable design for an axial rotary magnet assembly for magnetocaloric cooling. Magnetostatic simulations are performed in an isogeometric analysis framework for an exact and flexible representation of the geometry, where the rotor and stator are coupled with  harmonic mortaring. For efficient topology and shape optimization, analytical gradients have been used within the IGA framework.
To account for the influence of the MCM, simple bounds have been applied to estimate the permeability for (an)isotropic materials.

The main findings  can be summarized as follows:
\begin{itemize}
    \item Based on the topology optimization results, an inner rotary magnet assembly is proposed as  axial rotary design for commercial magnetocaloric refrigeration. To generate the large magnetic fields required for magnetocaloric cooling to date, it is suggested to build the geometry with five simple magnet segments and try to omit a continuous shaft. If necessary, this design can be extended to a co-rotational one.
    \item Inhomogeneous magnetic fields have shown to correlate with performance losses of the AMRs \citep{PerformanceMeasurementsOnALargeScale}. To avoid such problems, shape optimization has been applied to the rotor in order to homogenize the magnetic flux density in the AMRs. It has been shown that peaks in the magnetic flux density can be avoided for sufficiently small AMRs by adjusting the shape of the soft iron parts. 
    This is not the case for other prototypes given by the literature \citep{optimizedrefrigeration}.
    \item A comprehensive torque study has been performed for a parametrized stator with equally spaced AMRs. We have seen that a high and odd number of AMRs is in general favorable for lower torques. 
    \item Different positive and negative examples were evaluated with simplified assumptions for the MCM. Having a high number of AMRs leads to a sharper distinction between the (de)magnetization phases. The best results are obtained with anisotropic AMRs that enforce the magnetic flux in radial direction.
\end{itemize}

For future research, it is recommended to further increase the modeling effort of the AMRs with a multiscale model that considers measurement data of the MCM. Also, the mechanical and hydraulic design need to be worked out in more detail in order to enable the assembly. One important geometry parameter that has yet to be determined is the width of the air gap. For this, more information about performance, efficiency and cost are needed.

The system presented will serve as a foundation for a manufacturing design that will be assembled and tested by the authors.

\section*{Acknowledgement}

This work is supported by the joint DFG/FWF Collaborative Research Centre CREATOR (CRC -- TRR361/F90) at TU Darmstadt, TU Graz and JKU Linz as well as the Graduate School CE within the Centre for Computational Engineering at TU Darmstadt. MAGNOTHERM thanks the EIC-Accelerator MCD under grant agreement number 190183588 for their financial support. Many thanks also to Christian Vogel for the fruitful discussions.



\end{document}